\documentclass[reqno]{amsart}
\usepackage{amssymb}

\newtheorem{theorem}{Theorem}[section]
\newtheorem{proposition}[theorem]{Proposition}
\newtheorem{lemma}[theorem]{Lemma}
\newtheorem{corollary}[theorem]{Corollary}
\newtheorem*{lemma*}{Lemma}

\theoremstyle{definition}
\newtheorem*{definition}{Definition}

\theoremstyle{remark}
\newtheorem*{remark}{Remark}

\numberwithin{equation}{section}

\begin{document}

\title[Determinantal Construction of Polynomials]
{Determinantal Construction of Orthogonal Polynomials Associated
with Root Systems}

\author{J.F. van Diejen}
\address{
Instituto de Matem\'atica y F\'{\i}sica, Universidad de Talca,
Casilla 747, Talca, Chile}

\author{L. Lapointe}
\address{
Department of Mathematics and Statistics, McGill University,
Montr\'eal, Qu\'ebec H3A 2K6, Canada}

\author{J. Morse}
\address{
Department of Mathematics, University of Pennsylvania, 209 South
33rd Street, Philadelphia, PA 19104, USA}

\thanks{Work supported in part by the Fondo Nacional de Desarrollo
Cient\'{\i}fico y Tecnol\'ogico (FONDECYT) Grants \# 1010217 and
\# 7010217, the Programa Formas Cuadr\'aticas of the Universidad
de Talca, the C\'atedra Presidencial in Number Theory, and NSF
Grant \# 0100179.}

\date{March, 2002}

\begin{abstract}
We consider semisimple triangular operators acting in the
symmetric component of the group algebra over the weight lattice
of a root system. We present a determinantal formula for the
eigenbasis of such triangular operators. This determinantal
formula gives rise to an explicit construction of the Macdonald
polynomials and of the Heckman-Opdam generalized Jacobi
polynomials.
\end{abstract}

\maketitle \tableofcontents

\section{Introduction}\label{sec1}
The main objective of this work concerns the explicit computation
of families of orthogonal polynomials associated with root
systems. Key examples of the families under consideration are the
Macdonald polynomials
\cite{mac:symmetric1,mac:symmetric2,mac:orthogonal} and the
Heckman-Opdam generalized Jacobi polynomials
\cite{hec-sch:harmonic,opd:harmonic}. The origin of the
Heckman-Opdam polynomials lies in the harmonic analysis of simple
Lie groups, where they appear (for special parameter values) as
zonal spherical functions on compact symmetric spaces
\cite{hec-sch:harmonic,hel:geometric}. Other important
applications of these polynomials arise in mathematical physics,
where they are used to express the eigenfunctions of the quantum
Calogero-Sutherland one-dimensional many-body systems
\cite{sut:exact1,sut:exact2,ols-per:quantum,hec-sch:harmonic}. The
Macdonald polynomials have similar applications: they appear as
zonal spherical functions on compact quantum symmetric spaces
\cite{nou:macdonalds,sug:zonal}, and they are used to express the
eigenfunctions of Ruijsenaars' ($q$-)difference
Calogero-Sutherland systems
\cite{rui:complete,die:diagonalization}. Depending on the specific
application of interest, our work may thus be viewed as providing
an explicit construction for the zonal spherical functions on
compact (quantum) symmetric spaces or for the eigenfunctions of
the (difference) Calogero-Sutherland type quantum many-body
models.

The usual definition of the Heckman-Opdam and Macdonald
polynomials involves a Gram-Schmidt type orthogonalization of the
monomial basis with respect to a generalized Haar measure
\cite{hec-sch:harmonic,opd:harmonic,mac:symmetric2,mac:orthogonal}.
This definition, although most appropriate from a theoretical
point of view, is not very adequate for the explicit computation
of the polynomials in question. The main result of this paper is a
determinantal formula for the Heckman-Opdam and Macdonald
polynomials that gives rise to an efficient recursive procedure
from which their expansion in the monomial basis can be
constructed explicitly. For the type $A$ root systems the
Heckman-Opdam polynomials reduce (in essence) to Jack's
polynomials \cite{sta:some,mac:symmetric1} and the Macdonald
polynomials reduce to Macdonald's symmetric functions
\cite{mac:symmetric1}. In this case the determinantal construction
of the polynomials under consideration was laid out in previous
work by Lapointe, Lascoux, and Morse
\cite{lap-las-mor:determinantal1,lap-las-mor:determinantal2}. More
specifically, the results of the present paper constitute a
generalization of the methods of Refs.
\cite{lap-las-mor:determinantal1,lap-las-mor:determinantal2} to
the case of arbitrary root systems. For the Heckman-Opdam families
we consider general (not necessarily reduced) root systems and
general values of the root multiplicity parameters. For the
Macdonald families, however, we restrict for technical reasons to
those (reduced) root systems for which the dual root system
$R^\vee$ has a minuscule weight (thus including the types $A_N$,
$B_N$, $C_N$, $D_N$, $E_6$ and $E_7$ while excluding the types
$BC_n$, $E_8$, $F_4$ and $G_2$).

The paper is organized as follows. Section \ref{sec2} introduces
the concept of a triangular operator in the Weyl-group invariant
component of the group algebra over the weight lattice of a root
system. In Section \ref{sec3} we present a method for
diagonalizing such triangular operators by means of a
determinantal formula. The Heckman-Opdam and Macdonald polynomials
are defined in Section \ref{sec4}. We employ the determinantal
formula from Section \ref{sec3} to build explicit expressions for
the monomial expansions of these polynomials in Section \ref{sec5}
(Heckman-Opdam) and Sections \ref{sec6}, \ref{sec7} (Macdonald),
respectively. For completeness, some technicalities concerning the
explicit evaluation of the determinant of a Hessenberg matrix are
recalled in Appendix \ref{appA} at the end of the paper. To
facilitate explicit computations, we have furthermore included a
useful formula for the calculation of the orders of stabilizer
subgroups of the Weyl group in Appendix \ref{appB}.

\section{Triangular Operators in the Symmetrized Group Algebra}\label{sec2}
In this section we define the concept of a triangular operator in
the Weyl-group invariant component of the group algebra over the
weight lattice of a root system. For preliminaries on root systems
the reader is referred to e.g. Refs.
\cite{bou:groupes,hum:introduction}.

Let $E$, $\langle \cdot ,\cdot \rangle$ be a real Euclidean space
spanned by an irreducible root system $R$ with Weyl group $W$. We
write $\mathcal{Q}$ and $\mathcal{Q}^+$  for the root lattice and
its nonnegative semigroup generated by the positive roots $R^+$
\begin{equation}
\mathcal{Q}= \text{Span}_\mathbb{Z} (R),\;\;\;\mathcal{Q}^+=
\text{Span}_\mathbb{N} (R^+) .
\end{equation}
The weight lattice $\mathcal{P}$ and the cone of dominant weights
$\mathcal{P}^+$ are given by
\begin{equation}
\mathcal{P}= \{ \lambda\in E \mid  \langle \lambda ,\alpha^\vee
\rangle \in\mathbb{Z},\; \forall \alpha\in R \}
\end{equation}
and
\begin{equation}
\mathcal{P}^+= \{ \lambda\in E \mid  \langle \lambda ,\alpha^\vee
\rangle \in\mathbb{N},\; \forall \alpha\in R^+ \} ,
\end{equation}
where $\alpha^\vee=2\alpha /\langle \alpha ,\alpha \rangle$. The
weight lattice is endowed with the natural partial order
\begin{equation}\label{po}
\lambda \succeq \mu \;\Longleftrightarrow \;\lambda -\mu \in
\mathcal{Q}^+.
\end{equation}

Let $\mathcal{Q}^\vee$ denote the dual root lattice generated by
the dual root system $R^\vee =\{ \alpha^\vee \mid \alpha\in R\}$.
The group algebra over the weight lattice
$\mathbb{R}[\mathcal{P}]$ is the algebra generated by the formal
exponentials $e^\lambda$, $\lambda\in\mathcal{P}$ subject to the
multiplication relation $e^\lambda e^\mu = e^{\lambda +\mu}$. This
algebra can be realized explicitly as the algebra $\mathcal{A}$ of
(Fourier) polynomials on the torus $\mathbb{T}=E/(2\pi
\mathcal{Q}^\vee)$ through the identification
\begin{equation}
e^\lambda = e^{i\langle \lambda , x\rangle},\;\;\;\;\;\; \lambda
\in\mathcal{P}
\end{equation}
(with $x\in \mathbb{T}$). Symmetrization with respect to the
action of the Weyl group produces the basis of monomial symmetric
functions $\{ m_\lambda \}_{\lambda\in\mathcal{P}^+}$ for the
space $\mathcal{A}^W$ of Weyl-group invariant polynomials on
$\mathbb{T}$, where
\begin{equation}
m_\lambda = \sum_{\mu\in W(\lambda)} e^\mu ,\;\;\;\;\;\; \lambda
\in\mathcal{P},
\end{equation}
with $W(\lambda)$ denoting the orbit of $\lambda$ with respect to
the action of the Weyl group.

We write $\mathcal{A}^W_\lambda$ for the finite-dimensional
highest weight subspace of $\mathcal{A}^W$ with highest weight
$\lambda\in\mathcal{P}^+$, i.e.,
$\mathcal{A}^W_\lambda=\text{Span}\{
m_\mu\}_{\mu\in\mathcal{P}^+,\, \mu\preceq\lambda}$.
\begin{definition}
A linear operator $D:\mathcal{A}^W\to\mathcal{A}^W$ is called {\em
triangular} if $D(\mathcal{A}^W_\lambda)\subseteq
\mathcal{A}^W_\lambda$ for all $\lambda\in\mathcal{P}^+$.
\end{definition}

\section{Determinantal Diagonalization}\label{sec3}
The triangularity of a linear operator $D$ in $\mathcal{A}^W$
reduces its eigenvalue problem to a finite-dimensional one. In
this section we diagonalize the triangular operators through a
determinantal representation of the eigenfunctions.

Let $D$ be a triangular operator and let $\{ s_\lambda
\}_{\lambda\in\mathcal{P}^+}$ be a second basis of $\mathcal{A}^W$
that is related to the monomial basis by a unitriangular
transformation:
\begin{equation}\label{2ndbase}
m_\lambda = \sum_{\mu\in\mathcal{P}^+,\,\mu\preceq\lambda}
a_{\lambda \mu}\, s_\mu ,\;\;\;\;\; a_{\lambda \lambda }=1
\end{equation}
($\lambda\in\mathcal{P}^+$). The triangularity implies that the
expansion of $D m_\lambda$ in the basis $\{ s_\lambda
\}_{\lambda\in\mathcal{P}^+}$ is of the form
\begin{equation}\label{Daction}
D\, m_\lambda = \sum_{\mu\in\mathcal{P}^+,\,\mu\preceq\lambda}
b_{\lambda \mu}\, s_\mu ,\;\;\;\; b_{\lambda
\lambda}=\epsilon_\lambda ,
\end{equation}
with the diagonal matrix elements $\epsilon_\lambda$,
$\lambda\in\mathcal{P}^+$ being precisely the eigenvalues of $D$.
\begin{definition}
The triangular operator $D$ is called {\em regular} if
$\epsilon_\mu\neq\epsilon_\lambda$ when $\mu \prec \lambda$.
\end{definition}
For a regular triangular operator the eigenvalues
$\epsilon_\lambda$, $\lambda\in\mathcal{P}^+$ are semisimple. Let
$\{ p_\lambda \}_{\lambda\in\mathcal{P}^+}$ be a corresponding
basis of eigenfunctions diagonalizing $D$. Clearly, $p_\lambda$
has a monomial expansion of the form
\begin{equation}\label{mexp}
p_\lambda = \sum_{\mu\in\mathcal{P}^+,\,\mu\preceq\lambda}
c_{\lambda \mu}\, m_\mu ,\;\;\;\;\; c_{\lambda \lambda }=1,
\end{equation}
where we have normalized such that $p_\lambda$ is monic. The
following theorem provides an explicit determinantal formula for
$p_\lambda$, given the action of $D$ on $m_\lambda$  expressed in
the basis $s_\lambda$, i.e., given the expansion coefficients
$a_{\lambda \mu}$ and $b_{\lambda \mu}$ in Eqs. \eqref{2ndbase}
and \eqref{Daction}.

\begin{theorem}[Determinantal Formula]\label{df:thm}
Let $D$ be a regular triangular operator in $\mathcal{A}^W$ whose
action on the monomial symmetric functions is given by Eqs.
\eqref{2ndbase} and \eqref{Daction}. Then the monic basis $\{
p_\lambda\}_{\lambda\in\mathcal{P}^+}$ of $\mathcal{A}^W$
diagonalizing $D$, in the sense that
\begin{equation*}
D p_\lambda =\epsilon_\lambda\, p_\lambda ,\;\;\;\;\; \forall
\lambda\in\mathcal{P}^+,
\end{equation*}
is given explicitly by the (lower) Hessenberg determinant
  \begin{equation*}
  p_\lambda =\frac{1}{\mathcal{E}_\lambda}
  \begin{vmatrix}
  m_{\lambda^{(1)}} & \epsilon_{\lambda^{(1)}}-\epsilon_{\lambda^{(n)}}  & 0 & \hdots & \hdots &
0 \\
  m_{\lambda^{(2)}} &  d_{\lambda^{(2)}\lambda^{(1)}}& \epsilon_{\lambda^{(2)}}-\epsilon_{\lambda^{(n)}} &
   0 &\hdots & 0 \\
\vdots & \vdots & & \ddots & \ddots & \vdots \\
\vdots & \vdots & \vdots & & \ddots & 0 \\
\makebox[1ex]{} m_{\lambda^{(n-1)}} &
d_{\lambda^{(n-1)}\lambda^{(1)}}& d_{\lambda^{(n-1)}\lambda^{(2)}}
&\cdots & &
\epsilon_{\lambda^{(n-1)}}-\epsilon_{\lambda^{(n)}} \\
m_{\lambda^{(n)}} &  d_{\lambda^{(n)}\lambda^{(1)}}&
d_{\lambda^{(n)}\lambda^{(2)}} & \cdots
 &\cdots &
d_{\lambda^{(n)}\lambda^{(n-1)}}
  \end{vmatrix} .
  \end{equation*}
Here $\lambda^{(1)}<\lambda^{(2)}<\cdots
<\lambda^{(n-1)}<\lambda^{(n)}=\lambda$ denotes any linear
ordering of the dominant weights $\{ \mu\in\mathcal{P}^+ \mid
\mu\preceq\lambda \}$ refining the natural order \eqref{po}, the
normalization is determined by
\begin{equation*}
\mathcal{E}_\lambda =
\prod_{\mu\in\mathcal{P}^+,\,\mu\prec\lambda}
(\epsilon_\lambda-\epsilon_\mu) ,
\end{equation*}
and the matrix elements $d_{\lambda^{(j)}\lambda^{(k)}}$ ($n\geq
j>k\geq 1$) read
\begin{equation*}
d_{\lambda^{(j)}\lambda^{(k)}}=b_{\lambda^{(j)}\lambda^{(k)}}-
  \epsilon_{\lambda}\, a_{\lambda^{(j)}\lambda^{(k)}}.
\end{equation*}
\end{theorem}

\begin{proof}
Expansion of the determinant with respect to the first column
produces a linear combination of monomials in the highest weight
space $\mathcal{A}^W_\lambda$. The coefficient of the leading
monomial $m_\lambda$ is given by $(-1)^{n-1}$ times the product of
the elements on the super-diagonal, which are nonzero by the
regularity condition on $D$. Division by $\mathcal{E}_\lambda$
thus gives rise to a monic polynomial. It remains to show that
this polynomial is an eigenfunction of $D$ with eigenvalue
$\epsilon_\lambda$. To this end one observes that the action of
$(D-\epsilon_\lambda)$ on the determinant affects only its first
column. Indeed, we get---upon invoking the expansions
\eqref{2ndbase} and \eqref{Daction}---that
\begin{equation*}
  \bigl( D-\epsilon_{\lambda} \bigr)  p_{\lambda} =
  \frac{1}{\mathcal{E}_\lambda}
  \begin{vmatrix}
  & \ldots & \sum_{k=1}^{j-1}d_{\lambda^{(j)}\lambda^{(k)}}
  s_{\lambda^{(k)}} + (\epsilon_{\lambda^{(j)}}-
  \epsilon_{\lambda}) s_{\lambda^{(j)}} & \ldots & \\
  &&&&\\
  &\ldots& d_{\lambda^{(j)}\lambda^{(1)}}& \ldots&\\
  && \vdots & &\\
  &\ldots& d_{\lambda^{(j)}\lambda^{(j-1)}}& \ldots&\\
  &&&&\\
  &\ldots& \epsilon_{\lambda^{(j)}} -\epsilon_{\lambda} & \ldots &\\
  &&&&\\
  &\ldots& 0 & \ldots &\\
  && \vdots && \\
  \end{vmatrix}
\end{equation*}
(where, for typographical reasons, we have taken the transpose of
our matrix). The latter determinant has a first row of the form
$s_{\lambda^{(1)}}\text{row}_2+
s_{\lambda^{(2)}}\text{row}_3+\cdots+
s_{\lambda^{(n-1)}}\text{row}_{n}$, and thus vanishes identically.
\end{proof}

As a corollary of the determinantal formula for $p_\lambda$, one
arrives at a linear recurrence relation encoding an efficient
algorithm for the computation of the coefficients $c_{\lambda\mu}$
entering the monomial expansion \eqref{mexp}.

\begin{corollary}[Linear Recurrence Relation]\label{lrr:cor}
The monomial expansion of $p_\lambda$ is of the form
\begin{equation*}
p_\lambda = \sum_{\ell=1}^{n} c_{\lambda\lambda^{(\ell)}}\,
m_{\lambda^{(\ell)}},
\end{equation*}
with $c_{\lambda\lambda^{(n)}}=c_{\lambda \lambda}=1$ and
\begin{equation*}
 c_{\lambda\lambda^{(\ell-1)}}
=\frac{1}{\epsilon_\lambda-\epsilon_{\lambda^{(\ell-1)}}}
\sum_{k=\ell}^n c_{\lambda\lambda^{(k)}} \,
d_{\lambda^{(k)}\lambda^{(\ell-1)}}
\end{equation*}
($1<\ell \leq n$).
\end{corollary}
\begin{proof}
Immediate from Theorem \ref{df:thm} and the Hessenberg determinant
evaluation given by the lemma in Appendix \ref{appA} at the end of
the paper.
\end{proof}

Moreover, by solving the recurrence relation we arrive at the
following explicit expression for the coefficients
$c_{\lambda\mu}$ of the monomial expansion \eqref{mexp}.

\begin{corollary}[Explicit Monomial Expansion]\label{eme:cor}
The coefficients of the monomial expansion $p_\lambda =
\sum_{\ell=1}^{n} c_{\lambda\lambda^{(\ell)}}\,
m_{\lambda^{(\ell)}}$ are given explicitly by
\begin{equation*}
 c_{\lambda\lambda^{(\ell)}}
= \sum_{\begin{subarray}{c} \ell=j_r<j_{r-1}<\cdots <j_1<j_{0}= n\\
r= 1,\ldots ,n-\ell\end{subarray}} \frac{
d_{\lambda^{(j_0)}\lambda^{(j_1)}}
d_{\lambda^{(j_1)}\lambda^{(j_{2})}}\cdots
d_{\lambda^{(j_{r-1})}\lambda^{(j_r)}}}
{(\epsilon_\lambda-\epsilon_{\lambda^{(j_1)}})\cdots
(\epsilon_\lambda-\epsilon_{\lambda^{(j_r)}})} ,
\end{equation*}
with the convention that empty sums are equal to $1$ (so $
c_{\lambda\lambda^{(n)}}=c_{\lambda\lambda}=1$).
\end{corollary}

\begin{proof}
In view of Corollary \ref{lrr:cor}, it suffices to check that the
stated expression for $ c_{\lambda\lambda^{(\ell)}}$ represents
the (unique) solution to the linear recurrence relation of
Corollary \ref{lrr:cor}, subject to the initial condition $
c_{\lambda\lambda^{(n)}} =1$. Firstly, the convention that empty
sums are equal to $1$ guarantees that the initial condition is
satisfied. Secondly, by isolating the last factor in the numerator
and denominator, it is seen that for $1< \ell \leq n$
\begin{eqnarray*}
 c_{\lambda\lambda^{(\ell-1)}} &=&
\sum_{\begin{subarray}{c} \ell-1=j_r<j_{r-1}<\cdots <j_1<j_{0}= n\\
r= 1,\ldots ,n-\ell+1\end{subarray}}
\frac{d_{\lambda^{(j_0)}\lambda^{(j_1)}}
d_{\lambda^{(j_1)}\lambda^{(j_{2})}}\cdots
d_{\lambda^{(j_{r-1})}\lambda^{(j_{r})}}}
{(\epsilon_\lambda-\epsilon_{\lambda^{(j_1)}})\cdots
(\epsilon_\lambda-\epsilon_{\lambda^{(j_{r})}})}
\\
&=&
\frac{1}{\epsilon_\lambda-\epsilon_{\lambda^{(\ell-1)}}}\sum_{k=\ell}^n
 \:\Biggl( d_{\lambda^{(k)}\lambda^{(\ell-1)}}  \times \\
&& \makebox[2em]{} \sum_{\begin{subarray}{c}
k=j_{r-1}<j_{r-2}<\cdots <j_1< j_{0}=n\\ r-1= 1,\ldots
,n-k\end{subarray}} \frac{d_{\lambda^{(j_0)}\lambda^{(j_1)}}
d_{\lambda^{(j_1)}\lambda^{(j_{2})}}\cdots
d_{\lambda^{(j_{r-2})}\lambda^{(j_{r-1})}}}
{(\epsilon_\lambda-\epsilon_{\lambda^{(j_1)}})\cdots
(\epsilon_\lambda-\epsilon_{\lambda^{(j_{r-1})}})} \Biggr) \\
&=&
\frac{1}{\epsilon_\lambda-\epsilon_{\lambda^{(\ell-1)}}}\sum_{k=\ell}^n
 d_{\lambda^{(k)}\lambda^{(\ell-1)}} \,  c_{\lambda\lambda^{(k)}} .
\end{eqnarray*}
\end{proof}

\section{Orthogonal Polynomials}\label{sec4}
In this section the Heckman-Opdam and Macdonald polynomials are
defined.

Let $\Delta (x)$ be a positive continuous function on the torus
$\mathbb{T}/(2\pi \mathcal{Q}^\vee)$ that is invariant with
respect to the action of the Weyl group (i.e., $\Delta
(w(x))=\Delta (x)$ for all $w\in W$). We equip $\mathcal{A}^W$
with an inner product structure associated to the weight function
$\Delta$
\begin{equation}\label{ip}
\langle f , g\rangle_\Delta = \frac{1}{|W|} \int_\mathbb{T} f
\overline{g}\,\Delta\,\text{d} x \;\;\;\;\;\;\;
(f,g\in\mathcal{A}^W),
\end{equation}
where $\overline{g}$ denotes the complex conjugate of $g$ and
$|W|$ is the order of the Weyl group. Let $\{p_{\lambda,
\Delta}\}_{\lambda\in\mathcal{P}^+}$ be the basis of
$\mathcal{A}^W$ that is obtained from the monomial symmetric basis
$\{m_\lambda\}_{\lambda\in\mathcal{P}^+}$ through application of
the Gram-Schmidt process with respect to the partial order
$\succeq$ \eqref{po}. More specifically, by definition $p_{\lambda
,\Delta}$ is the polynomial of the form
\begin{subequations}
\begin{equation}
p_{\lambda,\Delta} = \sum_{\mu\in \mathcal{P}^+,\, \mu \preceq
\lambda} c_{\lambda ,\mu}(\Delta)\, m_\mu , \;\;\;\;\; c_{\lambda
,\lambda }(\Delta)=1,
\end{equation}
with coefficients $c_{\lambda ,\mu}(\Delta)$ such that
\begin{equation}
\langle p_{\lambda ,\Delta} ,m_\mu\rangle_\Delta
=0\;\;\;\;\;\text{for}\;\;\;\;\mu\in\mathcal{P}^+\;\;\text{with}\;\;
\mu\prec\lambda .
\end{equation}
\end{subequations}

In the general, the basis $\{ p_{\lambda,\Delta}
\}_{\lambda\in\mathcal{P}^+}$ is {\em not} orthogonal with respect
to the inner product $\langle \cdot ,\cdot\rangle_\Delta$
\eqref{ip}, as the natural order $\succeq$ \eqref{po} is not
linear (unless $\text{dim}(E)=1$). However, for two important
special choices of the weight function $\Delta$ it has been shown
that the above partial Gram-Schmidt process indeed {\em does}
produce an orthogonal basis
\cite{hec-sch:harmonic,opd:harmonic,mac:symmetric2,mac:orthogonal}:
\begin{subequations}
\begin{equation}
\Delta^{HO} = \prod_{\alpha\in R} (1-e^\alpha)^{g_\alpha},
\label{how}
\end{equation}
with $g_\alpha \geq 0$ such that $g_{w(\alpha)}=g_\alpha$,
$\forall w\in W$, and
\begin{equation}
 \Delta^{M} =
\prod_{\alpha\in R} \frac{(e^\alpha;q)_\infty}{(t_\alpha
e^\alpha;q)_\infty} , \label{mw}
\end{equation}
\end{subequations}
with  $(z;q)_\infty =\prod_{m=0}^\infty (1-zq^m)$ and
$0<q,t_\alpha<1$ such that $t_{w(\alpha )}=t_\alpha$, $\forall
w\in W$.
\begin{definition}
The orthogonal polynomials $p_{\lambda,\Delta}$,
$\lambda\in\mathcal{P}^+$ associated to the weight functions
$\Delta^{HO}$ \eqref{how} and $\Delta^{M}$ \eqref{mw} are called
the {\em Heckman-Opdam polynomials}
\cite{hec-sch:harmonic,opd:harmonic} and the {\em Macdonald
polynomials} \cite{mac:symmetric2,mac:orthogonal}, respectively.
\end{definition}

\begin{remark}[i] In the above definition of the Heckman-Opdam
polynomials we may allow for a root system $R$ that is {\em
nonreduced}. For the Macdonald polynomials, on the other hand, we
always assume that $R$ be {\em reduced}.
\end{remark}

\begin{remark}[ii]
The weight functions $\Delta^{HO}$ and $\Delta^M$ are invariant
with respect to the action of the Weyl group on the variable (i.e.
$\Delta (w(x))=\Delta (x)$) because of the $W$-invariance of the
orbit parameters $g_\alpha$ and $t_\alpha$. The $W$-invariance of
these parameters moreover implies that the values of $g_\alpha$
and $t_\alpha$ depend only on the length of the root $\alpha$.
\end{remark}

\section{Heckman-Opdam Polynomials}\label{sec5}
In this section we apply the formalism of Section \ref{sec3} to
arrive at a determinantal construction of the Heckman-Opdam
polynomials for arbitrary (not necessarily reduced) root systems.

\subsection{The hypergeometric differential operator}
To a vector $x\in E$ we associate the directional derivation
$\partial_x$ in $\mathcal{A}$, whose action on the exponential
basis is given by
\begin{equation}
\partial_x e^\lambda = \langle\lambda ,x  \rangle\, e^\lambda
\;\;\;\;\;\; (\lambda\in\mathcal{P}).
\end{equation}
\begin{definition}[\cite{hec-sch:harmonic,opd:harmonic}]
Let $x_1,\ldots ,x_N$ be an orthonormal basis of $E$. The
second-order partial differential operator
\begin{equation}\label{hdo}
  D = \sum_{j=1}^N \partial_{x_j}^2  +
  \sum_{\alpha\in R^+} g_\alpha
  \left( \frac{1+e^{-\alpha}}{1-e^{-\alpha}} \right)
  \partial_\alpha
\end{equation}
is called the {\em hypergeometric differential operator}
associated to the root system $R$.
\end{definition}
Clearly the definition of $D$ \eqref{hdo} does not depend on the
particular choice for the orthonormal basis $x_1,\ldots ,x_N$. It
is known that the hypergeometric differential operator maps the
space of invariants $\mathcal{A}^W$ into itself and, furthermore,
that the Heckman-Opdam polynomials form a basis of eigenfunctions
on which the operator acts diagonally
\cite{hec-sch:harmonic,opd:harmonic}.

We will now compute the action of $D$ \eqref{hdo} on the basis of
monomial symmetric functions. To this end some notation is needed.
We denote by $W_\lambda\subset W$ the stabilizer subgroup of
$\lambda\in\mathcal{P}$
\begin{equation*}
W_\lambda = \{ w\in W \mid w(\lambda)=\lambda \} ,
\end{equation*}
by $\rho_g\in E$ the weighted half-sum of the positive roots
\begin{equation*}
\rho_g =\frac{1}{2}\sum_{\alpha\in R^+} g_\alpha \,\alpha ,
\end{equation*}
by $r_\alpha:E\to E$ the orthogonal reflection in the hyperplane
perpendicular to $\alpha\in R$ through the origin
\begin{equation*}
r_\alpha (x)=x-\langle x, \alpha^\vee \rangle \alpha \;\;\;\;\;\;
(x\in E) ,
\end{equation*}
and by $[ p]$ the function that extracts the integral part of a
nonnegative real number $p$ through truncation.

\begin{lemma}[Action of the Hypergeometric Differential Operator]
\label{hdo-action:lem} The action of $D$ \eqref{hdo} on
$m_\lambda$, $\lambda\in\mathcal{P}^+$, is given by
\begin{eqnarray*}
  D\, m_\lambda &=& \bigl( \langle \lambda +\rho_g ,\lambda +\rho_g  \rangle
  -\langle \rho_g ,\rho_g \rangle  \bigr)\, m_\lambda   \\
 && \;\;\; + \;\frac{1}{|W_\lambda|} \sum_{\alpha\in R^+}
   \Bigl( g_\alpha\, \langle \lambda ,\alpha\rangle
  \sum_{\ell=1}^{[\langle \lambda ,\alpha^\vee\rangle/2]}
   |W_{\lambda -\ell\alpha}| \, |W^\alpha (\lambda -\ell\alpha )|\,
m_{\lambda -\ell\alpha} \Bigr) , \nonumber
\end{eqnarray*}
where $W^\alpha \subset W $ denotes the subgroup of order $2$
generated by $r_\alpha $ (so $| W^\alpha  (\lambda-\ell\alpha ) |$
is equal to $1$ if $\ell=\langle \lambda ,\alpha^\vee\rangle /2$
and equal to $2$ otherwise).
\end{lemma}

\begin{proof}
The computation of the expansion of $D m_\lambda$ in the monomial
basis hinges on the fundamental identity
\begin{eqnarray}
 \lefteqn{\left( \frac{1+e^{-\alpha}}{1-e^{-\alpha}} \right)
  \partial_\alpha
\bigl( e^\lambda + e^{r_\alpha (\lambda)}
\bigr) } \nonumber \\
&& = \langle \lambda ,\alpha \rangle \, e^\lambda\,
(1+e^{-\alpha}) \Bigl( \frac{ 1 - e^{-\langle \lambda, \alpha^\vee
\rangle \, \alpha}}{1-e^{-\alpha}}
\Bigr) \nonumber \\
&& = \langle \lambda ,\alpha \rangle \bigl( e^\lambda +
e^{\lambda-\langle \lambda, \alpha^\vee \rangle \, \alpha } \bigr
)   +\; 2\, \langle \lambda ,\alpha \rangle \sum_{\ell=1}^{\langle
\lambda ,\alpha^\vee\rangle -1}
e^{\lambda -\ell\alpha} \nonumber \\
&& = \langle \lambda ,\alpha \rangle \bigl( e^\lambda +
e^{r_\alpha (\lambda)} \bigr ) \nonumber\\ && \;\; \;\;\; + \;
\langle \lambda ,\alpha \rangle \sum_{\ell=1}^{[\langle \lambda
,\alpha^\vee\rangle /2]} |W^\alpha (\lambda -\ell\alpha)|\, \bigl(
e^{\lambda -\ell\alpha}+e^{r_\alpha(\lambda -\ell\alpha)} \bigr)
.\label{id}
\end{eqnarray}
Indeed, the following sequence of elementary manipulations reduces
the computation of the action of the first-order component of $D$
\eqref{hdo} on $m_\lambda$ to an application of identity
\eqref{id}:
\begin{eqnarray*}
\lefteqn{ \sum_{\alpha\in R^+} g_\alpha
  \left( \frac{1+e^{-\alpha}}{1-e^{-\alpha}} \right)
  \partial_\alpha m_\lambda } && \nonumber \\
&=& \frac{1}{|W_\lambda|} \sum_{\alpha\in R^+} g_\alpha
  \Bigl( \frac{1+e^{-\alpha}}{1-e^{-\alpha}} \Bigr)
  \partial_\alpha
\sum_{w\in W} e^{w(\lambda)} \nonumber\\
&=&\frac{1}{|W_\lambda|} \sum_{w\in W} w \Bigl( \sum_{\alpha\in
R^+} g_\alpha
  \Bigl( \frac{1+e^{-\alpha}}{1-e^{-\alpha}} \Bigr)
  \partial_\alpha
 e^{\lambda}\Bigr) \nonumber \\
&=&\frac{1}{2 |W_\lambda|} \sum_{w\in W} w \Bigl( \sum_{\alpha\in
R^+} g_\alpha
  \Bigl( \frac{1+e^{-\alpha}}{1-e^{-\alpha}} \Bigr)
  \partial_\alpha
\bigl( e^{\lambda}+ e^{r_\alpha(\lambda)}\bigr) \Bigr)\nonumber \\
&\stackrel{\text{Eq.~\eqref{id}}}{=}&\frac{1}{|W_\lambda|}
\sum_{\begin{subarray}{c} w\in W \\ \alpha \in R^+\end{subarray}}
w \Bigl( g_\alpha \langle \lambda ,\alpha \rangle \, e^{\lambda }
 + \, g_\alpha \, \langle \lambda ,\alpha \rangle
\sum_{\ell=1}^{[\langle \lambda ,\alpha^\vee\rangle/2]}
|W^\alpha (\lambda -\ell\alpha )|\, e^{\lambda -\ell\alpha}\Bigr)\nonumber \\
&=&
 \sum_{\alpha\in R^+}  g_\alpha \langle \lambda ,\alpha \rangle
\Biggl(  m_\lambda + \,\frac{1}{|W_\lambda|}
  \sum_{\ell=1}^{[\langle \lambda ,\alpha^\vee\rangle /2]} \,
|W_{\lambda -\ell\alpha}|\, |W^\alpha (\lambda -\ell\alpha )| \,
m_{\lambda -\ell\alpha} \Biggr) .
\end{eqnarray*}
Combined with the action of the second-order component of $D$
\eqref{hdo} on $m_\lambda$
\begin{equation*}
\sum_{j=1}^N \partial_{x_j}^2\, m_\lambda = \langle \lambda
,\lambda \rangle \, m_\lambda ,
\end{equation*}
this produces the formula of the lemma.
\end{proof}

It is a standard property of root systems that for any $\lambda\in
\mathcal{P}^+$ the integral convex hull $\mathcal{P}_\lambda= \{
\mu\in\mathcal{P} \mid W(\mu ) \preceq \lambda \}$ is saturated,
i.e., if $\mu\in\mathcal{P}_\lambda$ then $\mu-\ell\alpha\in
\mathcal{P}_\lambda$ for every integer $\ell$ between $0$ and
$\langle \mu ,\alpha^\vee\rangle$ (extremal values included)
\cite{hum:introduction}. Hence, it follows from Lemma
\ref{hdo-action:lem} that the hypergeometric differential operator
is triangular. To compute for $\mu$ dominant the coefficient of
$m_\mu$ in $D\, m_\lambda$, it suffices to collect all terms in
the lemma  for which $\lambda -\ell\alpha \in W(\mu)$. Notice in
this connection that for a given $\alpha\in R^+$ the
$\alpha$-string $\lambda -\alpha,\lambda-2\alpha,\ldots, \lambda-
[\langle \lambda ,\alpha^\vee \rangle/2]\alpha$ may hit the Weyl
orbit of $\mu$ at most once. Indeed, it is clear from expanding
both sides of the equality $||\lambda
-\ell^\prime\alpha||^2=||\lambda -\ell \alpha ||^2$ that $\lambda
-\ell^\prime \alpha \in W(\lambda -\ell\alpha)$---with $1\leq
\ell^\prime ,\ell\leq [\langle \lambda
,\alpha^\vee\rangle/2]$---implies $\ell^\prime =\ell$.

We thus end up with the following explicit triangular matrix
representation of the hypergeometric differential operator with
respect to the monomial basis.

\begin{proposition}[Triangular Expansion]\label{tmr:prp}
Let $\lambda\in \mathcal{P}^+$. We have that
\begin{equation*}
D\, m_\lambda = \epsilon_\lambda\, m_\lambda +
\sum_{\mu\in\mathcal{P}^+,\,\mu\prec\lambda} b_{\lambda \mu}\,
m_\mu ,
\end{equation*}
with
\begin{eqnarray*}
\epsilon_\lambda &=& \langle \lambda +\rho_g ,\lambda +\rho_g
\rangle
  -\langle \rho_g ,\rho_g \rangle  , \\
b_{\lambda \mu} &=&
 \frac{|W_\mu|}{|W_\lambda|}
 \sum_{\alpha\in [\lambda ,\mu]}
g_\alpha \,\langle \lambda ,\alpha \rangle\, n_{\lambda
\mu}(\alpha) .
\end{eqnarray*}
Here $[\lambda ,\mu] $ denotes the subset of roots $\alpha \in
R^+$ for which $\lambda -\ell\alpha\in W(\mu)$ for some (unique)
$\ell\in \{ 1,2,\ldots , [\langle \lambda ,\alpha^\vee \rangle
/2]\}$, and
\begin{equation*}
n_{\lambda \mu}(\alpha )=
\begin{cases}
1 & if\;\;  ||\mu||=|| P_{\alpha}(\lambda)||  \\
2 & if\;\;  ||\mu|| \neq  || P_{\alpha}(\lambda)||
\end{cases} ,
\end{equation*}
where $P_\alpha=(Id+r_\alpha)/2$ is the orthogonal projection onto
the hyperplane perpendicular to $\alpha$ through the origin. (So
for $\lambda -\ell\alpha\in W(\mu)$ we have that $n_{\lambda
\mu}(\alpha )=1$ if $\ell=\langle \lambda ,\alpha^\vee \rangle /2$
and $n_{\lambda \mu}(\alpha )=2$ otherwise.)
\end{proposition}

The regularity of the hypergeometric differential operator is
immediate from the following proposition.
\begin{proposition}[Monotonicity]\label{mon:prp}
For nonnegative parameters $g_\alpha $, the eigenvalues
$\epsilon_\lambda = \langle \lambda +\rho_g ,\lambda +\rho_g
\rangle
  -\langle \rho_g ,\rho_g \rangle $ are strictly monotonous in
$\lambda\in\mathcal{P}^+$, i.e.,
\begin{equation*}
\forall \lambda ,\mu \in \mathcal{P}^+:\;\;\;\; \mu \prec \lambda
\Longrightarrow \epsilon_\mu < \epsilon_\lambda .
\end{equation*}
\end{proposition}
\begin{proof}
Assume $\lambda ,\mu$ dominant with $\mu \prec \lambda$, and let
$\nu=\lambda -\mu$ (so $\nu \in \mathcal{Q}^+$). Then
\begin{equation*}
\epsilon_\lambda -\epsilon_\mu = \langle \nu ,\nu \rangle +
2\langle \mu+\rho_g ,\nu\rangle ,
\end{equation*}
which is positive in view of the fact that $\langle \nu
,\nu\rangle >0$ and $\langle \mu+\rho_g ,\nu\rangle\geq 0$ (since
both $\mu$ and $\rho_g$ lie in the closure of the dominant Weyl
chamber, cf. remark below).
\end{proof}

\begin{remark}
In the proof of Proposition \ref{mon:prp}, we used the fact that
for nonnegative parameters $g_\alpha$ the weighted half-sum
$\rho_g$ lies in the closure of the dominant Weyl chamber $\{ x\in
E \mid \langle x,\alpha \rangle >0,\;\forall\alpha\in R^+\}$. This
follows from the fact that, for any root $\beta$, the partial
half-sum $\rho (\beta)$ of positive roots with length $||\beta ||$
lies in the dominant cone $\mathcal{P}^+$:
\begin{equation*}
\rho(\beta) = \frac{1}{2}\sum_{\begin{subarray}{c}
                  \alpha\in R^+ \\
                  ||\alpha||=||\beta||
                 \end{subarray}} \alpha
=\sum_{\begin{subarray}{c}
        j=1,\ldots ,N \\
        ||\alpha_j||=||\beta||
       \end{subarray}} \omega_j,
\end{equation*}
where $\{ \alpha_j\}_{1\leq j\leq N}$ denotes the basis of simple
roots generating $\mathcal{Q}^+$ and $\{\omega_j\}_{1\leq j\leq
N}$ is the corresponding dual basis of fundamental weights
generating $\mathcal{P}^+$, such that $\langle
\omega_j,\alpha_k^\vee\rangle =\delta_{j,k}$. Indeed, since any
simple reflection $r_{\alpha_j}$ permutes the positive roots other
than $\alpha_j$ \cite{hum:introduction}, one has that
\begin{equation*}
r_{\alpha_j}(\rho(\beta))=
\begin{cases}
\rho (\beta) -\alpha_j & \text{if}\;\; ||\alpha_j||=||\beta||, \\
\rho (\beta)           & \text{otherwise} .
\end{cases}
\end{equation*}
It thus follows, from working out both sides of the equality
$\langle r_{\alpha_j} (\rho(\beta)),
r_{\alpha_j}(\alpha_j^\vee)\rangle = \langle \rho (\beta ),
\alpha_j^\vee\rangle $, that $\langle \rho (\beta),
\alpha_j^\vee\rangle $ is equal to $1$ if $||\alpha_j||
=||\beta||$ and is equal to $0$ otherwise. This entails that
$\rho(\beta)=\sum_{j=1}^N \langle \rho(\beta),\alpha_j^\vee\rangle
\omega_j =\sum^N_{j=1,\,
||\alpha_j||=||\beta||}\omega_j\in\mathcal{P}^+$ as claimed.
\end{remark}

\subsection{Determinantal construction}\label{sub52}
It is known that the hypergeometric differential operator $D$
\eqref{hdo} is symmetric with respect to the inner product
$\langle\cdot ,\cdot\rangle_\Delta$ \eqref{ip}, associated to the
weight function $\Delta^{HO}$ \eqref{how}. Combined with the
triangularity, this implies that the eigenbasis $\{ p_\lambda
\}_{\lambda\in\mathcal{P}^+}$ diagonalizing $D$ is given by the
Heckman-Opdam polynomials $p_{\lambda \,\Delta}$,
$\lambda\in\mathcal{P}^+$, defined in Section \ref{sec4} through
the (partial) Gram-Schmidt process
\cite{hec-sch:harmonic,opd:harmonic}. Moreover, since $D$ is
regular (Proposition \ref{mon:prp}), and its triangular action on
the monomial basis is known explicitly (Proposition
\ref{tmr:prp}), we can in fact construct this eigenbasis in closed
form by means of the determinantal construction in Section
\ref{sec3} (with $s_\lambda =m_\lambda$, so $a_{\lambda \mu}=1$ if
$\mu=\lambda$ and $a_{\lambda \mu}=0$ otherwise). This gives rise
to the following explicit representation of the Heckman-Opdam
polynomials.

\begin{theorem}[Determinantal Construction]\label{ho-dc:thm}
For $\lambda\in\mathcal{P}^+$, let
\begin{equation*}
p_\lambda = m_\lambda +
\sum_{\mu\in\mathcal{P}^+,\,\mu\prec\lambda} c_{\lambda \mu}\,
m_\mu
\end{equation*}
denote the (monic) Heckman-Opdam polynomial with parameters
$g_\alpha \geq 0$. Then we have---upon setting for $\mu ,\nu\in
\mathcal{P}^+$
\begin{eqnarray*}
\epsilon_\mu &=& \langle \mu +\rho_g ,\mu +\rho_g  \rangle
  -\langle \rho_g ,\rho_g \rangle  , \\
d_{\mu \nu} &=&
 \frac{|W_\nu|}{|W_\mu|}
 \sum_{\alpha\in [\mu ,\nu]}
g_\alpha \,\langle \mu ,\alpha \rangle\, n_{\mu \nu}(\alpha) ,
\end{eqnarray*}
with $[\mu ,\nu] \subset R^+$ and $n_{\mu \nu}(\alpha )$ in
accordance with the definition in Proposition
\ref{tmr:prp}---that:
\begin{itemize}
\item[i)] the polynomial $p_\lambda$ is represented explicitly
by the determinantal formula in Theorem \ref{df:thm},
\item[ii)] the coefficients $c_{\lambda \mu}$ of its monomial
expansion are generated by the linear recurrence in Corollary
\ref{lrr:cor},
\item[iii)] the expansion coefficients $c_{\lambda \mu}$ are
given in closed form by the formula in Corollary \ref{eme:cor}.
\end{itemize}
\end{theorem}

Given a concrete root system $R$, Theorem \ref{ho-dc:thm} turns
into an efficient algorithm for the computation of the associated
Heckman-Opdam polynomials. We will illustrate this below for the
classical root systems.

\begin{remark}[i]
The orders of the stabilizers in Theorem \ref{ho-dc:thm} can be
computed by means of the formula
\begin{equation}
|W_\lambda | =
\prod_{\begin{subarray}{c} \alpha\in R^+ \\
             \langle \lambda ,\alpha^\vee\rangle = 0
       \end{subarray}}
\frac{\langle \rho ,\alpha^\vee\rangle
+1+\frac{1}{2}\delta_{\frac{\alpha}{2}}} {\langle \rho
,\alpha^\vee\rangle +\frac{1}{2}\delta_{\frac{\alpha}{2}} } ,
\end{equation}
where $\rho =\frac{1}{2}\sum_{\alpha \in R^+} \alpha$, and
$\delta_{\alpha}=1$ if $\alpha \in R$ and zero otherwise. This
expression can be found e.g. in Ref. \cite[Section
12]{mac:orthogonal}, where it appears as a special case of the
norm formulas for the Macdonald polynomials. For the reader's
convenience, we have included a short proof of this formula in
Appendix \ref{appB}.
\end{remark}

\begin{remark}[ii]
It is clear from Proposition \ref{tmr:prp} that the matrix
$b_{\lambda \mu}$, representing the action of the hypergeometric
differential operator with respect to the monomial basis, is quite
sparse. The same is therefore true for the matrix $d_{\mu \nu}$
appearing in the determinantal formula for the Heckman-Opdam
polynomials in Theorem \ref{ho-dc:thm}. This means in practice
that the algorithm for generating the Heckman-Opdam polynomials
with the aid of the determinantal construction turns out to be
much faster than one would expect based just on the size of the
matrices involved.
\end{remark}

\begin{remark}[iii]
It is immediate from Theorem \ref{ho-dc:thm} that the coefficients
in the monomial expansion of the monic Heckman-Opdam polynomial
$p_\lambda$ are of the form $c_{\lambda \mu}=p_{\lambda
\mu}(g_\alpha )/q_{\lambda \mu}(g_\alpha)$, where $p_{\lambda
\mu}(g_\alpha )$ and $q_{\lambda \mu}(g_\alpha)$ are polynomials
in the parameters $g_\alpha$ that have nonnegative integral
coefficients (and with the denominators $q_{\lambda
\mu}(g_\alpha)$ dividing the normalization factor
$\mathcal{E}_\lambda =
\prod_{\mu\in\mathcal{P}^+,\,\mu\prec\lambda} (\epsilon_\lambda
-\epsilon_\mu)$). Recently, a much stronger positive-integrality
result for these expansion coefficients was found by Sahi
\cite{sah:new} (see also \cite{mac:symmetric1,kno-sah:recursion}
for the case $R=A_N$).
\end{remark}

\begin{remark}[iv]
It is clear from the proof of Proposition \ref{mon:prp} that the
hypergeometric differential operator $D$ \eqref{hdo} is in fact
regular for generic (complex) parameters $g_\alpha$ such that
$\langle \nu ,\nu \rangle + 2\langle \mu+\rho_g ,\nu\rangle \neq
0$ for all $\nu\in\mathcal{Q}^+$ and $\mu\in\mathcal{P}^+$. Hence,
the determinantal construction of the Heckman-Opdam polynomials in
Theorem \ref{ho-dc:thm}, as the eigenbasis for the hypergeometric
differential operator, extends meromorphically to $g_\alpha$ in
the complex plane.
\end{remark}

\begin{remark}[v]
It is known that the coefficients of the Heckman-Opdam polynomials
can in principle be computed by means of (cumbersome) Freudenthal
type recurrence relations \cite{hec-sch:harmonic}. {}From this
perspective, the determinantal construction of Theorem
\ref{ho-dc:thm} thus provides the explicit solution to this
Freudenthal type recurrence. The recurrence in Part $ii)$ of
Theorem \ref{ho-dc:thm}---which arises as a particular case of the
general recurrence scheme in Corollary \ref{lrr:cor} upon choosing
for our triangular operator the hypergeometric differential
operator---reads concretely
\begin{eqnarray*}
\lefteqn{(\epsilon_{\lambda}-\epsilon_{\lambda^{(\ell-1)}}) \,
c_{\lambda \lambda^{(\ell-1)}}
= } && \\
&& \sum_{k= \ell}^n
\frac{|W_{\lambda^{(\ell-1)}}|}{|W_{\lambda^{(k)}}|}
 \sum_{\alpha\in [\lambda^{(k)} ,\, {\lambda^{(\ell-1)}}]}
g_\alpha \,\langle \lambda^{(k)} ,\alpha \rangle\,
n_{\lambda^{(k)} {\lambda^{(\ell-1)}}}(\alpha) \, c_{\lambda
\lambda^{(k)}} .
\end{eqnarray*}
It may in fact be seen as a suitable symmetric reduction of the
Freudenthal type recurrence relations, enabling their explicit
solution in closed form via Corollary \ref{eme:cor}. When
$g_\alpha =1$, $\forall \alpha\in R$, our recurrence is closely
related to the optimized Freudenthal recurrence scheme for the
computation of weight multiplicities of characters of simple Lie
groups due to Moody and Patera \cite{moo-pat:fast}.
\end{remark}

\subsection{Tables for the classical root systems}\label{HOtables}
We will now provide tables of the matrix elements building the
determinantal formulas for the Heckman-Opdam polynomials
associated with the classical root systems. In each case, we will
only list the minimum amount of information needed for
constructing the matrix, viz., {\em i.} the cone of the dominant
weights and its partial order, {\em ii.} the eigenvalues building
the super-diagonal of the matrix, {\em iii.} and the values of the
lower-triangular matrix elements. For further data on the root
systems of interest we refer to the tables in Bourbaki
\cite{bou:groupes}.

It will be convenient to parameterize the dominant weights of the
classical root systems in terms of $N$-tuples
\begin{subequations}
\begin{equation}
\lambda = (\lambda_1,\lambda_2, \ldots ,\lambda_N)
\end{equation}
of weakly decreasing (half-)integers (so $\lambda_1\geq
\lambda_2\geq \cdots \geq \lambda_N$). Often we think of these
$N$-tuples also as multi-sets of the form
\begin{equation}
\lambda = \{ \lambda_1 ,\lambda_2 ,\ldots ,\lambda_N \} ,
\end{equation}
\end{subequations}
where the parts $\lambda_j$ are listed from largest to smallest.
For $\lambda=(\lambda_1,\dots,\lambda_N)$ ($
=\{\lambda_1,\dots,\lambda_N\}$), we define $\lambda^\varepsilon=
(\lambda_1,\dots,\lambda_{N-1},\varepsilon |\lambda_N|)$ ($
=\{\lambda_1,\dots,\lambda_{N-1},\varepsilon |\lambda_N|\} $),
with $\varepsilon \in \{1,-1\}$. We need the following two
operations on our weakly decreasing $N$-tuples:
\begin{subequations}
\begin{eqnarray}
\lambda \setminus \mu &=&
 \{ \lambda_1 ,\lambda_2 ,\ldots ,\lambda_N \}\setminus
 \{ \mu_1 ,\mu_2 ,\ldots ,\mu_N \} , \\
\lambda \ominus \mu &=& ( \lambda^+\setminus \mu^+
,(\mu^+\setminus \lambda^+)^\varepsilon ),
\end{eqnarray}
\end{subequations}
where $\varepsilon = \text{sign}(\lambda_N)\times \text{sign}
(\mu_N)$. The first operation takes the difference of $\lambda$
and $\mu$ as multi-sets, i.e., taking into account the
multiplicities of the parts. (By convention, we will list the
parts of this difference again from large to small.) The second
operation encodes---up to a possible sign---the symmetric
difference of $\lambda^+$ and $\mu^+$. For example:
$(5,3,2\frac{1}{2},1,1)\ominus (4,3,3,1,-1)= (\{ 5,2\frac{1}{2} \}
, \{ 4,-3  \})$. Finally, for future reference we furthermore
introduce the operations
\begin{subequations}
\begin{eqnarray}
|\lambda | &=& \lambda_1+\lambda_2+\cdots +\lambda_N ,\\
\eta_{\lambda} (m) &=& | \{ j=1,\ldots ,N \mid \lambda_j=m\,
\vee\, \lambda_j=-m \} | ,\label{mult} \\
\bar{\lambda} &=&
(\lambda_1,\lambda_2,\ldots ,\lambda_{N-1},-\lambda_N)
,\label{conj}
\end{eqnarray}
\end{subequations}
producing, respectively, the sum of the parts, the number of parts
with specified absolute value $|m|$, and the conjugate $N$-tuple
with the sign of the last part flipped.

\subsubsection{The case $R=A_{N-1}$}
For the type $A$ root system the Heckman-Opdam polynomials reduce
(in essence) to Jack polynomials \cite{sta:some,mac:symmetric1}.
The determinantal construction in Theorem \ref{ho-dc:thm}
reproduces in this particular case the determinantal construction
of the Jack polynomials found by Lapointe, Lascoux, and Morse
\cite{lap-las-mor:determinantal2}.

In dealing with the type $A$ root system, it is more convenient to
work with partitions rather than with the weights themselves. Let
$\Lambda_N$ be the set of partitions with at most $N$ parts, i.e.,
the set of weakly decreasing $N$-tuples with components given by
nonnegative integers. For $\lambda ,\mu \in\Lambda_N$ the
dominance order on these partitions is defined as
\begin{equation}\label{no}
\lambda \succeq \mu \Longleftrightarrow |\lambda |
=|\mu|\;\;\text{and}\;\; \sum_{j=1}^k (\lambda_j-\mu_j) \geq 0\;
\text{for}\; k=1,\ldots ,N-1.
\end{equation}
We write $\hat{\lambda}$ for the orthogonal projection of
$\lambda\in\Lambda_N$ onto the hyperplane $E\subset \mathbb{R}^N$
perpendicular to the vector $(1,1,\ldots ,1)$:
\begin{equation}\label{proj}
\hat{\lambda}= (\lambda_1,\ldots ,\lambda_N) - \frac{|\lambda
|}{N} (1,\ldots ,1) .
\end{equation}
The cone of dominant weights associated to the root system
$A_{N-1}$ is now given by the projection of $\Lambda_N$ onto the
hyperplane $E$, i.e. $\mathcal{P}^+_A = \{ \hat{\lambda} \mid
\lambda\in\Lambda_N \} $, equipped with a partial order induced by
the dominance ordering of the partitions in Eq. \eqref{no}.
Specifically, for given $\lambda\in\Lambda_N$ all dominant weights
smaller or equal to $\hat{\lambda}\in\mathcal{P}^+_A$ are given by
\begin{equation}\label{poA}
\mathcal{P}^+_{\hat{\lambda},A}= \{ \hat{\mu}\mid
\mu\in\Lambda_N\; \wedge\; \mu\preceq \lambda \} .
\end{equation}
The projection $\lambda\to\hat{\lambda}$ \eqref{proj} has a
nontrivial kernel of the form $(1,1,\ldots ,1)\mathbb{N}$. The set
in Eq. \eqref{poA}, however, is clearly independent of the
particular choice for the partition $\lambda$ projecting onto the
dominant weight $\hat{\lambda}$.

The Weyl group acts transitively on the root system $A_{N-1}$, as
all roots have the same length. Thus, the value of the root
multiplicity parameter $g_\alpha$ does not depend on $\alpha$,
viz. $g_\alpha =g$ for all $ \alpha \in R$. Given a partition
$\lambda\in \Lambda_N$, let us define for $\mu\preceq\lambda$
\begin{subequations}
\begin{equation}
\epsilon_\mu^A = \sum_{j=1}^N \mu_j\, (\mu_j + g\,(N+1-2j)),
\end{equation}
and for $\nu\prec\mu\preceq\lambda$
\begin{equation}\label{dA}
d_{\mu\nu}^A =
\begin{cases}
2 g\, (m_1-m_2)\,\mathcal{N}_{\nu}(n_1,n_2) & \text{if}\;
\mu\ominus\nu= (\{ m_1 ,m_2\} ,\{n_1 ,n_2\} ) \\
& \text{with}\; m_1-n_1 =n_2-m_2 >0, \\
0 & \text{otherwise},
\end{cases}
\end{equation}
\end{subequations}
where
\begin{equation}\label{N}
\mathcal{N}_{\nu}(n_1 ,n_2) =
\begin{cases}
\eta_{\nu}(n_1) \eta_{\nu}(n_2) & \text{if}\; |n_1 |\neq |n_2|, \\
\binom{\eta_{\nu}(n_1)}{2} & \text{if}\; |n_1|=|n_2 |,
\end{cases}
\end{equation}
and $\eta_{\nu} (n)$ denotes the multiplicity counter defined in
Eq. \eqref{mult}. The super-diagonal
$\epsilon^A_{\hat{\mu}}-\epsilon^A_{\hat{\lambda}}$ ($\hat{\mu}
\preceq \hat{\lambda}$) and the lower triangular block
$d_{\hat{\mu}\hat{\nu}}^A$ ($\hat{\nu}\prec\hat{\mu} \preceq
\hat{\lambda}$) of the Hessenberg matrix in Theorem
\ref{ho-dc:thm} become for the $A_{N-1}$-type Heckman-Opdam
polynomial $p_{\hat{\lambda}}^A$:
\begin{equation}\label{projm}
\epsilon^A_{\hat{\mu}}-\epsilon^A_{\hat{\lambda}} =
\epsilon^A_{\mu}-\epsilon^A_{\lambda} \;\;\;\;\text{and}\;\;\;\;
d_{\hat{\mu}\hat{\nu}}^A= d_{\mu\nu}^A ,
\end{equation}
respectively. (Notice in this connection that the expressions
$\epsilon^A_{\mu}-\epsilon^A_{\lambda}$ and $d_{\mu\nu}^A$ on the
r.h.s. are invariant with respect to the additive action of
$(1,1,\ldots ,1)\mathbb{N}$ on $\Lambda_N$.)

\subsubsection{The case $D_{N}$}
The cone of dominant weights $\mathcal{P}^+_D$ consists of the
$N$-tuples $\lambda = (\lambda_1,\ldots ,\lambda_N)$ with parts
$\lambda_j$ that are all integers or all half-integers subject to
the ordering
\begin{equation}
\lambda_1\geq \lambda_2\geq \cdots \geq \lambda_{N-1}\geq
|\lambda_N| .
\end{equation}
The partial order on $\mathcal{P}_+^D$ is defined as
\begin{equation}
\lambda \succeq \mu \Longleftrightarrow
\begin{cases}
\sum_{j=1}^k (\lambda_j-\mu_j)\in \mathbb{N} &\text{for}\;
k=1,\ldots ,N-2, \\
\sum_{j=1}^{N-1} (\lambda_j-\mu_j) +\varepsilon
(\lambda_N-\mu_N)\in 2\mathbb{N} & \text{for}\; \varepsilon =\pm
1.
\end{cases}
\end{equation}

The Weyl group again acts transitively on the root system $D_{N}$,
so we have $g_\alpha= g$, $\forall\alpha\in R$. The super-diagonal
$\epsilon^D_{\mu}-\epsilon^D_{\lambda}$ ($\mu \preceq \lambda$)
and the lower triangular block $d_{\mu\nu}^D$ ($\nu\prec\mu
\preceq \lambda$) of the Hessenberg matrix in Theorem
\ref{ho-dc:thm} become for the $D_{N}$-type Heckman-Opdam
polynomial $p_{\lambda}^D$ of the form
\begin{subequations}
\begin{equation}
\epsilon_\mu^D = \sum_{j=1}^N \mu_j\, (\mu_j +2g\, (N-j))
\end{equation}
and
\begin{equation}
d_{\mu\nu}^D = \left \{
\begin{array}{l}
\left( d_{m_1,m_2;n_1,n_2}^A +d_{m_1,\overline{m_2};n_1,n_2}^A +
d_{m_1,m_2;n_1,\overline{n_2}}^A
+d_{m_1,\overline{m_2};n_1,\overline{n_2}}^A \right)
\mathcal{N}_{\nu}
(n_1,n_2) \\
\\
\qquad \qquad \qquad  \text{if}\; {\mu}\ominus {\nu}= (\{ m_1
,m_2\} ,\{n_1 ,n_2\} ) \text{~with~} n_{\mu}(0) \neq 0
\text{~and~} m_2 \neq 0, \\
\\
\left( d_{m_1,m_2;n_1,n_2}^A
+d_{m_1,\overline{m_2};n_1,\overline{n_2}}^A \right)
\mathcal{N}_{\nu}
(n_1,n_2) \\
\\
\qquad \qquad \qquad  \text{if}\; {\mu}\ominus {\nu}= (\{ m_1
,m_2\} ,\{n_1 ,n_2\} ) \text{~with~} n_{\mu}(0) = 0
\text{~or~} m_2 = 0, \\
\\
 d_{m,\overline{\Delta^+};\Delta^+,\overline{n}}^A \,\, \mathcal{N}_{\nu}
(\Delta^+,n) +  d_{m,\overline{\Delta^-};\Delta^-,n}^A \, \,
\mathcal{N}_{\nu}
(\Delta^-,n)\\
\\
\qquad \qquad \qquad  \text{if}\;
{\mu}\ominus {\nu}= (\{ m\} ,\{n\} ) \text{~with~} n_{\mu}(0) \neq 0,  \\
\\
 d_{m,\overline{\Delta^+};\Delta^+,\overline{n}}^A \,\, \mathcal{N}_{\nu}
(\Delta^+,n)
 \\
\\
\qquad \qquad \qquad  \text{if}\;
{\mu}\ominus {\nu}= (\{ m\} ,\{n\} ) \text{~with~} n_{\mu}(0) = 0 , \\
\\
0 \qquad \quad \qquad  \, \, \, \text{otherwise}.
\end{array} \label{dD}
\right.
\end{equation}
\end{subequations}
Here $d^A_{m_1,m_2;n_1,n_2}$ refers to the $A_1$-type matrix
elements (cf. Eq. \eqref{dA}), viz.,
\begin{equation}
d^A_{m_1,m_2;n_1,n_2} =
\begin{cases}
2\, g\, (m_1-m_2) & \text{if~} m_1-n_1=n_2-m_2 > 0,\\
0 & \text{otherwise},
\end{cases}
\end{equation}
and $\mathcal{N}_\nu(n_1,n_2)$ is the same as above (cf. Eq.
\eqref{N}). Furthermore,  $\overline{m}$ stands for $-m$ and
$\Delta^{\pm}= (m \pm n)/2$.

\begin{remark} In the first line of $d_{\mu \nu}^D$, at most two terms can
be nonzero if $n_2=0$, and at most one term otherwise. Similarly,
in the second line, at most one term can be nonzero.
\end{remark}

\subsubsection{The case $B_{N}$}
The cone of dominant weights $\mathcal{P}^+_B$ consists of the
$N$-tuples $\lambda = (\lambda_1,\ldots ,\lambda_N)$ with parts
$\lambda_j$ that are all integers or all half-integers subject to
the ordering
\begin{equation}
\lambda_1\geq \lambda_2\geq \cdots \geq \lambda_{N-1}\geq
\lambda_N\geq 0 .
\end{equation}
The partial order on $\mathcal{P}^+_B$ is defined as
\begin{equation}\label{BNpo}
\lambda \succeq \mu \Longleftrightarrow \sum_{j=1}^k
(\lambda_j-\mu_j)\in \mathbb{N}  \;\; \text{for}\; k=1,\ldots ,N.
\end{equation}

The $B_N$-type root system has two root lengths, so the action of
the Weyl group splits up in two orbits. We will denote the root
multiplicity parameters for the long and short roots by $g$ and
$g_s$, respectively. The super-diagonal
$\epsilon^B_{\mu}-\epsilon^B_{\lambda}$ ($\mu \preceq \lambda$)
and the lower triangular block $d_{\mu\nu}^B$ ($\nu\prec\mu
\preceq \lambda$) of the Hessenberg matrix in Theorem
\ref{ho-dc:thm} become for the $B_{N}$-type Heckman-Opdam
polynomial $p_{\lambda}^B$ of the form
\begin{subequations}
\begin{equation}
\epsilon_\mu^B = \sum_{j=1}^N \mu_j\, (\mu_j +2g\, (N-j)+g_s)
\end{equation}
and
\begin{equation}
d^B_{\mu\nu} =
\begin{cases}
d_{\mu\nu}^D + d_{\mu\nu}^{\text{short}} &\text{if}\; \mu =\bar{\mu} , \\
d_{\mu\nu}^D + d_{\bar{\mu}\nu}^D + d_{\mu\nu}^{\text{short}}
&\text{if}\; \mu \neq \bar{\mu}  ,
\end{cases}
\end{equation}
\end{subequations}
with $d_{\mu\nu}^D$ taken from Eq. \eqref{dD} and
\begin{equation}\label{ds}
d_{\mu\nu}^{\text{short}} =
\begin{cases}
2 g_s m\, \eta_{\nu}(n) & \text{if}\; \mu\ominus \nu = (\{ m\} ,\{n\} )  \\
& \text{with}\; m-n>0, \\
0& \text{otherwise} .
\end{cases}
\end{equation}

\subsubsection{The cases $C_N$ and $BC_{N}$}
The cone of dominant weights $\mathcal{P}^+_{BC}$ consists of the
partitions $\lambda = (\lambda_1,\ldots ,\lambda_N)$ in
$\Lambda_N$ (cf. the $A_{N-1}$-type above). The partial order on
$\mathcal{P}^+_{BC}$ is the same as in the $B_N$-case (cf. Eq.
\eqref{BNpo})
\begin{equation}\label{BCNpo}
\lambda \succeq \mu \Longleftrightarrow \sum_{j=1}^k
(\lambda_j-\mu_j)\geq 0  \;\; \text{for}\; k=1,\ldots ,N.
\end{equation}

The $BC_N$-type root system has three root lengths, so the action
of the Weyl group splits up into three orbits. We will denote the
root multiplicity parameters for the long and short roots by $g_l$
and $g_s$, respectively. The parameter for the remaining (i.e.
middle) roots is $g$. The super-diagonal
$\epsilon^{BC}_{\mu}-\epsilon^{BC}_{\lambda}$ ($\mu \preceq
\lambda$) and the lower triangular block $d_{\mu\nu}^{BC}$
($\nu\prec\mu \preceq \lambda$) of the Hessenberg matrix in
Theorem \ref{ho-dc:thm} become for the $BC_{N}$-type Heckman-Opdam
polynomial $p_{\lambda}^{BC}$ of the form
\begin{subequations}
\begin{equation}
\epsilon_\mu^{BC} = \sum_{j=1}^N \mu_j\, (\mu_j +2g\,
(N-j)+g_s+2\,g_l)
\end{equation}
and
\begin{equation}
d^{BC}_{\mu\nu} =
\begin{cases}
d_{\mu\nu}^D +
d_{\mu\nu}^{\text{short}} + d_{\mu\nu}^{\text{long}} &\text{if}\; \mu =\bar{\mu} ,\\
d_{\mu\nu}^D + d_{\bar{\mu}\nu}^D + d_{\mu\nu}^{\text{short}} +
d_{\mu\nu}^{\text{long}} &\text{if}\; \mu \neq \bar{\mu} ,
\end{cases}
\end{equation}
\end{subequations}
where $d_{\mu\nu}^D$ and $d_{\mu\nu}^{\text{short}}$ are taken
from Eqs. \eqref{dD} and \eqref{ds}, respectively, and
\begin{equation}
d_{\mu\nu}^{\text{long}} =
\begin{cases}
4 g_l m\, \eta_{\nu}(n) & \text{if}\; \mu\ominus \nu = (\{ m\} ,\{n\} ) \\
& \text{with}\; m-n\in 2\mathbb{N}, \\
0& \text{otherwise} .
\end{cases}
\end{equation}

\begin{remark}[i]
The $C_N$ case is obtained from the $BC_N$ case by setting
$g_s=0$. In this situation one generally can reduce the size of
the Hessenberg matrix, as the partial order on the weights for the
$C_N$ root system, viz.
\begin{equation}\label{CNpo}
\lambda\succeq \mu \Longleftrightarrow
\begin{cases}
\sum_{j=1}^k (\lambda_j-\mu_j)\in \mathbb{N} &\text{for}\;
k=1,\ldots ,N-1, \\
\sum_{j=1}^{N} (\lambda_j-\mu_j)  \in 2\mathbb{N} ,&
\end{cases}
\end{equation}
is less refined than the partial order in Eq. \eqref{BNpo}
corresponding to the $BC_N$ root system. More specifically, if the
monomial on the $l^{th}$ row is not comparable to the leading
monomial in the $C_N$ ordering \eqref{CNpo}, then we may eliminate
(for $g_s=0$) the $l^{th}$ row together with the $(l+1)^{th}$
column from the Hessenberg matrix. (To keep the normalization
monic, we should of course also delete the corresponding factors
from the normalization constant $\mathcal{E}_\lambda$.)
\end{remark}

\begin{remark}[ii]
The $BC_N$ Heckman-Opdam polynomials deserve special attention as
they are universal in the sense that the polynomials associated
with the other classical root systems can be obtained as special
cases (the types $B_N$, $C_N$ and $D_N$ by specialization of the
parameters $g_l$ and $g_s$, and the type $A_{N-1}$ by selecting
the top-degree homogeneous component). For a systematic study of
the properties of the $BC_N$-type Heckman-Opdam polynomials we
refer to Refs. \cite{bee-opd:certain,die:properties} and papers
cited therein.
\end{remark}

\begin{remark}[iii]
Example: for $R=B_3$ and $\lambda=(2,1,0)$ the determinantal
formula reads
\begin{footnotesize}
\begin{eqnarray*}
&& \makebox{\begin{normalsize}$p_{2,1,0}$\end{normalsize}} =
\bigl((2+4\,g)(1+2\,g+g_s)(3+4\,g+g_s)(4+6\,g+2\,g_s)(5+10\,g+3\,g_s)\bigr)^{-1}
\times\\
&& \begin{vmatrix}\makebox[0.5ex]{} m_{{0,0,0}}&
-5-10\,g-3\,g_{{s}} & 0 & 0 & 0 & 0 \\
\makebox[0.5ex]{}m_{{1,0,0}}&6\,g_{{s}}& -4- 6\,g-2\,g_{{s}}&
0 & 0 & 0 \\
\makebox[0.5ex]{}m_{{1,1,0}}&24\,g& 4\,g_{{s}}&-3-4\,g-g_{{s}}& 0 & 0 \\
\makebox[0.5ex]{}m_{{2,0,0}}& 12\,g_{{s}} & 4\,g_{{s}} &4\,g & -1-2\,g-g_{{s}} & 0 \\
\makebox[0.5ex]{}m_{{1,1,1}}&0 & 8\,g& 2\,g_{{s}}& 0& -2-4\,g\makebox[0.5ex]{}\\
\makebox[0.5ex]{}m_{{2,1,0}}& 0 & 24\,g+8\, g_{{s}}
&8\,g_{{s}}&4\,g_ {{s}}&12\,g
\end{vmatrix}.
\end{eqnarray*}
\end{footnotesize}
This polynomial may also be interpreted as a special case of the
$BC_3$-type Heckman-Opdam polynomial $p_{2,1,0}$ with $g_l=0$. We
observe in this connection that for $g_s=0$, the $1^{st}$, the
$3^{rd}$ and the $4^{th}$ row, together with the $2^{nd}$, the
$4^{th}$ and the $5^{th}$ column, may be eliminated from the
Hessenberg matrix (cf. Remark (i) above). Indeed, the weights
$(0,0,0)$, $(1,1,0)$ and $(2,0,0)$ are not comparable to the
highest weight $(2,1,0)$ with respect to the $C_N$-type partial
order in Eq. \eqref{CNpo}. (To keep our normalization monic, we
must also delete the $2^{nd}$, the $3^{rd}$, and the $5^{th}$
factor from the normalization constant.)
\end{remark}

\section{Macdonald Polynomials: the case $t_\alpha =t$}\label{sec6}
In this section we apply the formalism of Section \ref{sec3} to
arrive at a determinantal construction of the Macdonald
polynomials. Throughout this section it will be assumed that the
root system $R$ is reduced and that the dual root system $R^\vee$
has a minuscule weight (thus including the types $A_N$, $B_N$,
$C_N$, $D_N$, $E_6$ and $E_7$ while excluding the types $BC_n$,
$E_8$, $F_4$ and $G_2$). We will furthermore restrict to the case
that $t_\alpha =t$, $\forall \alpha\in R$ and---unless explicitly
stated otherwise---we will consider the $(q,t)$ parameters as
indeterminates rather than real (or complex) numbers.

\subsection{The Macdonald operator}
For $x\in E$, we define the $q$-translation in $\mathcal{A}$ via
its action on the exponential basis:
\begin{equation}
T_{x,q} \, e^\lambda = q^{\langle \lambda ,x\rangle}\, e^\lambda
\;\;\;\;\;\; (\lambda\in\mathcal{P}).
\end{equation}

\begin{definition}[\cite{mac:symmetric2,mac:orthogonal}]
Let $\pi$ be a minuscule weight for $R^\vee$, i.e., the vector
$\pi\in E$ is such that $\langle \pi ,\alpha \rangle \in \{ 0
,1\}$ for all $\alpha\in R^+$. The $q$-difference operator
\begin{equation}
D_\pi = \frac{1}{|W_\pi|} \sum_{w\in W} \Bigl( \prod_{\alpha\in
R^+} \frac{1-t^{\langle \pi ,\alpha\rangle} e^{w(\alpha)}}
     {1-e^{w(\alpha)}} \Bigr) T_{w(\pi), q}
\end{equation}
is called the {\em Macdonald operator} associated to the minuscule
weight $\pi$.
\end{definition}
(The above definition of the Macdonald operator $D_\pi$ is not
precisely the same as the one employed by Macdonald
\cite{mac:symmetric2,mac:orthogonal}; both definitions do coincide
upon restriction to the space of invariant polynomials
$\mathcal{A}^W$ though.) In order to compute the action of $D_\pi$
on the monomial basis we will make use of the Weyl characters
$\chi_\lambda$, $\lambda \in \mathcal{P}$:
\begin{subequations}
\begin{equation}
\chi_\lambda = \delta^{-1} \sum_{w\in W} \det (w)\, e^{w(\lambda
+\rho )} ,
\end{equation}
where $\rho$ and $\delta$ denote the half sum of the positive
roots and the Weyl denominator respectively
\begin{eqnarray}
\rho &=& \frac{1}{2} \sum_{\alpha\in R^+} \alpha , \\
\delta &=& \prod_{\alpha\in R^+} (e^{\alpha/2}-e^{-\alpha /2})
=\sum_{w\in W} \det (w)\, e^{w(\rho )} .
\end{eqnarray}
\end{subequations}
(Clearly the determinant $\det (w)$ is equal to $(-1)^{\ell(w)}$,
where $\ell(w)$ represents the length of the (shortest)
decomposition of $w$ into a product of simple reflections.) It is
well-known that for $\lambda\in \mathcal{P}^+$ one has that
\begin{equation}\label{kostka}
\chi_\lambda = \sum_{\mu\in\mathcal{P}^+,\,\mu\preceq\lambda}
K_{\lambda \mu}\, m_\mu ,\;\;\;\;\; K_{\lambda \lambda }=1,
\end{equation}
with coefficients $K_{\lambda \mu}\in \mathbb{N}$. (In fact, the
coefficients  $K_{\lambda \mu}$, which are also known as Kostka
numbers, count the multiplicity of the weight $\mu$ in the
irreducible representation of the Lie algebra corresponding to the
root system $R$ with highest weight $\lambda$.) An efficient way
to compute the coefficients $K_{\lambda \mu}$ is through the
application of Theorem \ref{ho-dc:thm} with $g_\alpha =1$,
$\forall \alpha\in R$. However, for our purposes such a
calculation is not necessary as we need the inverse of this basis
transformation rather than Eq. \eqref{kostka} itself (cf.
Corollary \ref{invkostka:cor} below).

It is evident from the expansion in Eq. \eqref{kostka} that the
Weyl characters $\{ \chi_\lambda \}_{\lambda\in\mathcal{P}^+}$
form a basis of $\mathcal{A}^W$. The following lemma provides a
formula for the action of the Macdonald operator $D_\pi$ on the
monomials $m_\lambda$ in terms of Weyl characters. For the root
system $R=A_N$ the formula in question is due to Macdonald
\cite{mac:symmetric1,mac:symmetric2}.

\begin{lemma}[Action of the Macdonald Operator]\label{Dpi-action:lem}
Let $\lambda\in\mathcal{P}^+$. Then one has that
\begin{equation*}
D_\pi\, m_\lambda =
 t^{\langle \pi ,\rho\rangle}
\sum_{\nu\in W(\lambda)} \Bigl( \sum_{\tau\in W(\pi)} t^{\langle
\tau, \rho \rangle} q^{\langle \tau, \nu \rangle} \Bigr) \chi_{\nu
} .
\end{equation*}
\end{lemma}

\begin{proof}
Our starting point is the Weyl denominator formula in the form
\begin{equation*}
e^{-\rho}\prod_{\alpha\in R^+} (e^\alpha -1)= \sum_{w\in W} \det
(w)\, e^{w(\rho )} .
\end{equation*}
By acting on both sides with the $t$-translator $T_{\pi ,t}$ we
obtain
\begin{equation*}
t^{-\langle \pi ,\rho\rangle} e^{-\rho}\prod_{\alpha\in R^+}
(t^{\langle \pi ,\alpha \rangle }e^\alpha -1)= \sum_{w\in W} \det
(w)\, t^{\langle \pi , w(\rho)\rangle} e^{w(\rho )} .
\end{equation*}
Division of the latter identity by the former gives rise to the
following expansion for the coefficients of the Macdonald operator
\begin{equation*}
\prod_{\alpha\in R^+} \frac{1-t^{\langle \pi ,\alpha\rangle}
e^{\alpha}}
     {1-e^{\alpha}} = \delta^{-1}
t^{\langle \pi ,\rho\rangle} \sum_{w\in W} \det (w)\, t^{\langle
\pi , w(\rho)\rangle} e^{w(\rho )} .
\end{equation*}
Substitution of this expansion in the definition of $D_\pi$
(taking into account the anti-symmetry of the Weyl denominator
$w(\delta)=\det(w)\,\delta$), and acting on the exponential
$e^\nu$ yields:
\begin{eqnarray*}
D_\pi \, e^{\nu} &=& \delta^{-1} \frac{t^{\langle \pi
,\rho\rangle}}{|W_\pi|} \sum_{w_1,w_2\in W} \det (w_1 w_2)\,
t^{\langle \pi ,w_2 (\rho)\rangle}
q^{\langle w_1(\pi),\nu\rangle} e^{\nu +w_1 w_2 (\rho)} \\
 &=&
 \delta^{-1} \frac{t^{\langle \pi ,\rho\rangle}}{|W_\pi|}
\sum_{w_1,w_2\in W} \det (w_1 w_2)\, t^{\langle w_1(\pi ) ,w_1 w_2
(\rho)\rangle}
q^{\langle w_1(\pi),\nu\rangle} e^{\nu +w_1 w_2 (\rho)} \\
 &\stackrel{w_1w_2\to w}{=}&
\delta^{-1} t^{\langle \pi ,\rho\rangle} \sum_{w\in W,\, \tau\in W
(\pi)} \det (w)\,
 t^{\langle \tau ,w (\rho)\rangle}
q^{\langle \tau , \nu\rangle} e^{\nu +w (\rho)} .
\end{eqnarray*}
Summation over $\nu$ in the orbit $W(\lambda)$ then entails:
\begin{eqnarray*}
D_\pi\, m_\lambda &=& \delta^{-1} t^{\langle \pi ,\rho\rangle}
\sum_{\begin{subarray}{c} w\in W \\
     \tau\in W (\pi),\, \nu\in W(\lambda)
      \end{subarray}}
\det (w)\,
 t^{\langle \tau ,w (\rho)\rangle}
q^{\langle \tau , \nu\rangle} e^{\nu +w (\rho)} \\
 &=&
\delta^{-1} t^{\langle \pi ,\rho\rangle}
\sum_{\begin{subarray}{c}w\in W \\
      \tau\in W (\pi),\, \nu\in W(\lambda)
      \end{subarray}}
\det (w)\,
 t^{\langle \tau ,w (\rho)\rangle}
q^{\langle \tau , w(\nu )\rangle} e^{w(\nu +\rho)} \\
&=& t^{\langle \pi ,\rho\rangle} \sum_{\begin{subarray}{c}
       \tau\in W(\pi) \\
       \nu\in W(\lambda)
       \end{subarray}}
t^{\langle \tau , \rho\rangle} q^{\langle \tau , \nu\rangle}
\chi_\nu ,
\end{eqnarray*}
which completes the proof.
\end{proof}

For $\lambda\in\mathcal{P}$, let $w_\lambda$ be the unique
shortest Weyl group element such that $w_\lambda (\lambda)\in
\mathcal{P}^+$. Then it follows from the definition of the Weyl
characters that for $\nu\in\mathcal{P}$
\begin{equation}
\chi_\nu =
\begin{cases}
\det (w_{\nu+\rho})\, \chi_{(w_{\nu+\rho}(\nu+\rho)-\rho)} &
\text{if}\; |W_{\nu+\rho}|=1 , \\
0 & \text{if}\; |W_{\nu+\rho}|> 1 .
\end{cases}
\end{equation}
(Notice in this connection that---in view of Corollary
\ref{stabilizer:cor} in Appendix \ref{appB}---the stabilizer of a
weight is nontrivial if and only if there exist a root $\alpha\in
R^+$ perpendicular to it, i.e., if and only if there exists a
reflection $r_\alpha$, $\alpha\in R^+$ stabilizing the weight in
question.) Hence, to find for $\lambda ,\mu\in\mathcal{P}^+$ the
multiplicity of $\chi_\mu$ in $D_\pi\, m_\lambda$, we have to
collect all terms in the formula of Lemma \ref{Dpi-action:lem}
corresponding to weights $\nu\in W(\lambda)$ such that
$w_{\nu+\rho }(\nu+\rho )-\rho =\mu$, or equivalently, $\nu=w_{\nu
+\rho}^{-1}(\mu +\rho)-\rho$. Clearly the action of $D_\pi$ is
triangular \cite{mac:symmetric2,mac:orthogonal}, since
\begin{equation*}
\mu = w_{\nu+\rho }(\nu+\rho )-\rho= w_{\nu+\rho
}(\nu)-(\rho-w_{\nu+\rho }(\rho))\preceq w_{\nu+\rho }(\nu)\preceq
w_\nu (\nu)=\lambda
\end{equation*}
(where in the two last steps we used the fact that any dominant
weight $\lambda$ is maximal in its Weyl orbit, i.e.,
$w(\lambda)\preceq \lambda$ for all $w\in W$
\cite{hum:introduction}). We thus arrive at the following explicit
triangular expansion of $D_\pi m_\lambda$ in terms of $\chi_\mu$.

\begin{proposition}[Triangular Expansion]\label{qtmr:prp}
Let $\lambda\in \mathcal{P}^+$. We have that
\begin{equation*}
D_\pi\, m_\lambda = \epsilon_\lambda\, \chi_\lambda +
\sum_{\mu\in\mathcal{P}^+,\,\mu\prec\lambda} b_{\lambda \mu}\,
\chi_\mu ,
\end{equation*}
with
\begin{eqnarray*}
\epsilon_\lambda &=&  t^{\langle \pi ,\rho\rangle} \sum_{\tau\in
W(\pi)} t^{\langle \tau, \rho \rangle}
q^{\langle \tau,\lambda \rangle} , \\
b_{\lambda \mu} &=& \sum_{\nu\in W(\lambda) \cap (W(\mu +\rho
)-\rho)} \det (w_{\rho +\nu})\, \epsilon_\nu  .
\end{eqnarray*}
\end{proposition}

For $t=1$ the Macdonald operator $D_\pi$ trivializes to
$\sum_{\tau\in W(\pi)} T_{\tau,\, q}$, which acts diagonally on
$m_\lambda$ through multiplication by the eigenvalue
$\sum_{\tau\in W(\pi)} q^{\langle \tau,\lambda \rangle}$. The
formula of Lemma \ref{Dpi-action:lem} reduces in this case (and
upon division by the eigenvalue) to the following well-known
relation between the symmetric monomials and the Weyl characters:
\begin{equation}
m_\lambda = \sum_{\mu\in W(\lambda)} \chi_{\mu } .
\end{equation}
In the same way, one recovers from Proposition \ref{qtmr:prp} the
inverse of the expansion in Eq. \eqref{kostka}.
\begin{corollary}[Inverse Kostka Numbers]\label{invkostka:cor}
Let $\lambda\in\mathcal{P}^+$. The expansion of the symmetric
monomial $m_\lambda$ in terms of Weyl characters is given by
\begin{equation*}
m_\lambda = \chi_\lambda +
\sum_{\mu\in\mathcal{P}^+,\,\mu\prec\lambda} a_{\lambda \mu}\,
\chi_\mu ,
\end{equation*}
with
\begin{equation*}
a_{\lambda \mu} =\sum_{\nu\in W(\lambda) \cap (W(\mu +\rho
)-\rho)} \det (w_{\rho +\nu}) .
\end{equation*}
\end{corollary}

The following proposition guarantees that the Macdonald operator
$D_\pi$ is regular.
\begin{proposition}[Regularity]\label{qreg:prp}
The Macdonald operator $D_\pi$ is regular in the sense that
\begin{equation*}
\forall \lambda ,\mu \in \mathcal{P}^+:\;\;\;\; \mu \prec \lambda
\Longrightarrow \epsilon_\mu (q,t) \neq \epsilon_\lambda (q,t)
\end{equation*}
(as (analytic) functions of the indeterminates $q$ and $t$).
\end{proposition}
\begin{proof}
After setting $t=q^g$ and $q=\exp (z)$, we get
\begin{eqnarray*}
t^{-\langle \pi ,\rho\rangle}\epsilon_\lambda  &=& \sum_{\tau\in
W(\pi)}
\exp (z\langle \tau ,\lambda +g\rho\rangle ) \\
&=& |W(\pi)|+ z\sum_{\tau\in W(\pi)} \langle \tau ,\lambda
+g\rho\rangle + \frac{z^2}{2} \sum_{\tau\in W(\pi)} \langle \tau
,\lambda +g\rho\rangle^2
+ O(z^3) \\
&=& |W(\pi)|+ c_\pi z^2 \langle \lambda +g\rho ,\lambda
+g\rho\rangle + O(z^3) ,
\end{eqnarray*}
with $c_\pi >0$. (In the last step we employed the fact that the
$W$-invariant linear form $\sum_{\tau\in W(\pi)} \langle \tau ,x
\rangle $ vanishes and that the $W$-invariant positive quadratic
form $\sum_{\tau\in W(\pi)} \langle \tau ,x \rangle^2$ must be
proportional to $\langle x,x\rangle $, because the representation
of the Weyl group on $E$ is irreducible and unitary). When $g$ is
positive, one has that $\langle \mu +g\rho ,\mu +g\rho\rangle <
\langle \lambda +g\rho ,\lambda +g\rho\rangle$ for all dominant
weights $\mu ,\lambda$ with $\mu \prec \lambda$ in view of
Proposition \ref{mon:prp}. It thus follows that for comparable
dominant weights $\mu\neq\lambda$ the corresponding eigenvalues
$\epsilon_\mu (q,t)$ and $\epsilon_\lambda (q,t) $ cannot be equal
as (analytic) functions of the indeterminates $q$ and $t$.
\end{proof}

\subsection{Determinantal construction}
We will now apply the determinantal formalism of Section
\ref{sec3} to construct the eigenbasis of $D_\pi$. To this end we
pick for the second basis $\{ s_\lambda\}_{\lambda\in
\mathcal{P}^+}$ the basis of Weyl characters $\{
\chi_\lambda\}_{\lambda\in \mathcal{P}^+}$. Specifically, by
plugging in the eigenvalues $\epsilon_\lambda$ and off-diagonal
matrix elements $b_{\lambda \mu}$ from Proposition \ref{qtmr:prp},
together with the inverse Kostka numbers $a_{\lambda \mu}$ from
Corollary \ref{invkostka:cor}, the formulas of Theorem
\ref{df:thm} and the Corollaries \ref{lrr:cor} and \ref{eme:cor}
give rise to the desired eigenbasis of the corresponding Macdonald
operator $D_\pi$. For parameters such that $0<q,t<1$, this
eigenbasis coincides with the Macdonald polynomials defined in
Section \ref{sec4} through the (partial) Gram-Schmidt process
\cite{mac:symmetric2,mac:orthogonal}. We thus end up with the
following determinantal construction of the Macdonald polynomials
(not necessarily with $0<q,t<1$).

\begin{theorem}[Determinantal Construction]\label{m-dc:thm}
For $\lambda\in\mathcal{P}^+$, let
\begin{equation*}
p_\lambda = m_\lambda +
\sum_{\mu\in\mathcal{P}^+,\,\mu\prec\lambda} c_{\lambda \mu}\,
m_\mu
\end{equation*}
denote the (monic) Macdonald polynomial with $t_\alpha=t$,
$\forall \alpha\in R$. Then we have---upon setting for $\mu
,\nu\in\mathcal{P}^+$
\begin{eqnarray*}
\epsilon_\mu &=& t^{\langle \pi ,\rho\rangle} \sum_{\tau\in
W(\pi)} t^{\langle \tau, \rho \rangle}
q^{\langle \tau,\mu \rangle}  , \\
d_{\mu \nu} &=&  \sum_{\kappa\in W(\mu) \cap (W(\nu +\rho )-\rho)}
\det (w_{\rho +\kappa})\, (\epsilon_\kappa -\epsilon_\lambda )
\end{eqnarray*}
(so $d_{\mu\mu}=\epsilon_\mu-\epsilon_\lambda$)---that:
\begin{itemize}
\item[i)] the polynomial $p_\lambda$ is represented explicitly
by the determinantal formula in Theorem \ref{df:thm},
\item[ii)] the coefficients $c_{\lambda \mu}$ of its monomial
expansion are generated by the linear recurrence in Corollary
\ref{lrr:cor},
\item[iii)] the expansion coefficients $c_{\lambda \mu}$ are
given in closed form by the formula in Corollary \ref{eme:cor}.
\end{itemize}
\end{theorem}

\begin{remark}[i]
To determine the matrix elements $d_{\mu \nu}$, it is not
efficient to compute the intersection $W(\mu)\cap
(W(\nu+\rho)-\rho)$ for each $\mu ,\nu \in\mathcal{P}^+$ such that
$\nu\prec \mu\preceq\lambda$. Indeed, because the matrices at
issue are sparse, a better strategy is to construct the matrix row
by row. For this purpose one first determines for each dominant
weight $\mu\preceq\lambda$ the set $\Lambda_\mu= \{ \tilde{\kappa}
\in\mathcal{P} \mid \tilde{\kappa}\in \rho+W(\mu),\;
|W_{\tilde{\kappa}}|=1\}$ (i.e. all regular points of the
translated Weyl orbit $\rho+W(\mu)$). Weyl-permuting the weights
in $\Lambda_\mu$ to the dominant cone and translating over $-\rho$
produces all the nonzero contributions to the row $\mu$.
Specifically, the nonzero matrix elements on the row corresponding
to $\mu$ occur in the columns corresponding to $\nu$ from the set
$\{ \nu\in\mathcal{P}^+ \mid \nu = w_{\tilde{\kappa}}
(\tilde{\kappa} )-\rho ,\; \tilde{\kappa}\in \Lambda_\mu \}$. The
matrix elements in question are built  of contributions of the
form $\det(w_{\tilde{\kappa}}) (\epsilon_{\tilde{\kappa}
-\rho}-\epsilon_\lambda )$, $\tilde{\kappa}\in\Lambda_\mu$.
\end{remark}

\begin{remark}[ii]
It is well-known that for $t=q^g$ and $q\to 1$ the Macdonald
polynomial $p_\lambda$ tends to the corresponding Heckman-Opdam
polynomial (with $g_\alpha =g$, $\forall \alpha\in R$)
\cite{mac:orthogonal}. To perform this limit at the level of the
above determinantal construction, it suffices to determine the
asymptotics of the eigenvalues $\epsilon_\mu$ for $q\to 1$. The
asymptotics in question is given by (cf. the proof of Proposition
\ref{qreg:prp})
\begin{equation}\label{qasymptotics}
\epsilon_\mu t^{-\langle \pi ,\rho\rangle} = |W(\pi)| + c_\pi
\langle \mu +g\rho,\mu +g\rho\rangle (q-1)^2+ O((q-1)^3),
\end{equation}
where $c_\pi $ is a positive constant that does not depend on
$\mu$ and $g$. Since the formulas of the determinantal
construction for the Macdonald polynomials are invariant with
respect to an affine rescaling of the spectrum of the form
$\epsilon_\mu\to a\epsilon_\mu +b$ (with $a\neq 0$), we only pick
up the second-order term of the asymptotics in Eq.
\eqref{qasymptotics} when sending $q$ to $1$. The upshot is that
by replacing $\epsilon_\mu$ by $\langle \mu +g\rho,\mu
+g\rho\rangle$ in Theorem \ref{m-dc:thm}, we wind up with an
alternative determinantal formula for the Heckman-Opdam
polynomials (with $g_\alpha =g$, $\forall \alpha\in R$, and with
$R$ such that $R^\vee$ has a minuscule weight). From a practical
point of view the formulas coming from Theorem \ref{m-dc:thm} are
less efficient than those of Theorem \ref{ho-dc:thm}, however, as
the action of the Macdonald operator expanded in Weyl characters
tends to be much less sparse than the action of the hypergeometric
differential operator expanded in monomials. As a consequence, the
matrices entering the determinantal formulas for the Heckman-Opdam
polynomials coming from Theorem \ref{m-dc:thm} are much less
sparse than those of the type given by Theorem \ref{ho-dc:thm}.
\end{remark}

\begin{remark}[iii]
The recurrence relation in Part $ii)$ of Theorem \ref{m-dc:thm}
reads concretely
\begin{eqnarray*}
\lefteqn{(\epsilon_{\lambda}-\epsilon_{\lambda^{(\ell-1)}}) \,
c_{\lambda \lambda^{(\ell-1)}} =} && \\
&& \sum_{k= \ell}^n  \, \, \sum_{\kappa\in W({\lambda^{(k)}}) \cap
(W({\lambda^{(\ell-1)}}
 +\rho )-\rho)}
\det (w_{\rho +\kappa})\, (\epsilon_\kappa -\epsilon_\lambda ) \,
c_{\lambda \lambda^{(k)}} .
\end{eqnarray*}
This relation should be regarded as a symmetrized Freudenthal type
recurrence for the coefficients in the monomial expansion of the
Macdonald polynomials. For $\epsilon_\mu = \langle \mu + g\rho
,\mu +g\rho\rangle$, this recurrence degenerates to a recurrence
for the coefficients of the Heckman-Opdam polynomials with
$g_\alpha =g$, $\forall \alpha \in R$ (cf. Remark (ii) above). The
recurrence in question is different from the previous recurrence
for the Heckman-Opdam polynomials originating from the
hypergeometric differential operator (cf. Remark (v) at the end of
Section \ref{sub52}). In particular, for $g=1$ this gives rise to
an alternative system of symmetrized Freudenthal type recurrence
relations for the weight multiplicities of characters of simple
Lie groups.
\end{remark}

\subsection{Tables for the classical root systems}\label{Mtables}
We now provide tables for the construction of the Macdonald
polynomials associated with the classical root systems. The
minimum information needed for invoking Theorem \ref{m-dc:thm}
consists of {\em i.} the cone of dominant weights and its partial
order, {\em ii.} the half-sum of the positive roots $\rho$, {\em
iii.} the action of the Weyl group, and {\em iv.} the eigenvalues
$\epsilon_\mu$. Below we list items {\em ii.}--{\em iv} for
$R=A_{N-1}$, $B_N$, $C_N$ and $D_N$. For item {\em i} the reader
is referred to Subsection \ref{HOtables}.

\subsubsection{The case $A_{N-1}$}
The $A_{N-1}$-type Macdonald polynomials amount (in essence) to
the Macdonald symmetric functions of Ref. \cite{mac:symmetric1}.
Theorem \ref{m-dc:thm} reproduces in this case the determinantal
construction of the Macdonald symmetric functions due to Lapointe,
Lascoux, and Morse \cite{lap-las-mor:determinantal1}.

We will again formulate the construction in terms of partitions
with at most $N$ parts, by adding a trivial center to the weight
lattice (cf. Subsection \ref{HOtables}). The Weyl group $W$ is
given by the permutation group of $N$ letters $\Sigma_N$. A Weyl
group element $w=\sigma\in \Sigma_N$ acts on a partition
$\lambda\in\Lambda_N$ by rearranging its parts
\begin{equation}
\sigma (\lambda_1,\ldots ,\lambda_N) = (\lambda_{\sigma_1},\ldots
,\lambda_{\sigma_N}) .
\end{equation}

To construct the Macdonald symmetric function associated to a
partition $\lambda\in\Lambda_N$, one employs Theorem
\ref{m-dc:thm} with $\mathcal{P}^+=\Lambda_N$ endowed with the
dominance order in Eq. \eqref{no}, and
\begin{subequations}
\begin{eqnarray}
\epsilon_\mu^A &=& \sum_{j=1}^N t^{N-j} q^{\mu_j} , \\
 \rho_A &=& (N-1,N-2,\ldots ,1,0) .
\end{eqnarray}
\end{subequations}
For a rearrangement $\kappa$ of a partition $\mu\in\Lambda_N$, the
sign $\det (w_{\rho+\kappa})$ is given by the signature of the
shortest permutation $\sigma_{\rho+\kappa }$ rearranging
$\rho+\kappa$ such that its parts become weakly decreasing (i.e.
$\det (w_{\rho+\kappa})=(-1)^{\ell(\sigma_{\rho+\kappa})}$, where
$\ell(\sigma_{\rho+\kappa})$ denotes the number of transpositions
of the permutation).

Projection of the resulting Macdonald symmetric function onto the
space of homogeneous functions of degree zero (i.e., replacing
$m_\mu$ by $m_{\hat{\mu}}$ in the monomial expansion) entails the
$A_{N-1}$-type Macdonald polynomial $p^A_{\hat{\lambda}}$
associated to the weight $\hat{\lambda}$ \eqref{proj}.

\subsubsection{The case $B_N$}
The Weyl group is the semi-direct product of $\Sigma_N$ and the
$N$-fold product of $\mathbb{Z}_2=\mathbb{Z}/(2\mathbb{Z})$, i.e.,
$W=\Sigma_N\ltimes \mathbb{Z}_2^N$. A Weyl group element
$w=(\sigma ,\varepsilon)$ acts on a weight
$\lambda\in\mathcal{P}_B^+$ as
\begin{equation}\label{waction}
w(\lambda_1,\ldots ,\lambda_N) =
(\varepsilon_1\lambda_{\sigma_1},\ldots
,\varepsilon_N\lambda_{\sigma_N}),
\end{equation}
with $\varepsilon_j\in \{ 1 ,-1\}$ for $j=1,\ldots ,N$.

To construct the $B_N$-type Macdonald polynomial $p_\lambda^{B}$
associated to a weight $\lambda\in\mathcal{P}^B_+$, one employs
Theorem \ref{m-dc:thm} with
\begin{subequations}
\begin{eqnarray}
\epsilon_\mu^B &=& \sum_{j=1}^N \left (t^{2N-j} q^{\mu_j}
+ t^{j-1} q^{-\mu_j} \right) ,\\
 \rho_B &=&
(N-\frac{1}{2},N-\frac{3}{2},\ldots ,\frac{3}{2},\frac{1}{2}) .
\end{eqnarray}
\end{subequations}
For a weight $\kappa\in W(\mu)$ with $\mu\in\mathcal{P}^+_B$, the
sign $\det (w_{\rho+\kappa})$ is given by the signature of the
shortest permutation rearranging $\rho+\kappa$ such that the
absolute values of its parts become weakly decreasing, multiplied
by $(-1)^{n_\varepsilon}$ where
\begin{equation}
n_\varepsilon = |\{ j=1,\ldots , N \mid \rho_j+\kappa_j < 0 \}| .
\end{equation}

\subsubsection{The case $C_N$}
The Weyl group and its action on a weight
$\lambda\in\mathcal{P}^+_C$ are the same as in the $B_N$-case.

To construct the $C_N$-type Macdonald polynomial $p_\lambda^{C}$
associated to a weight $\lambda\in\mathcal{P}^C_+$, one employs
Theorem \ref{m-dc:thm} with
\begin{subequations}
\begin{eqnarray}
\epsilon^C_\mu &=& \prod_{j=1}^N \left(
t^{N+1-j}q^{\mu_j/2}+q^{-\mu_j/2}\right) ,\\
\rho_C &=&  (N,N-1,\ldots ,2,1) .
\end{eqnarray}
\end{subequations}
The sign $\det (w_{\rho+\kappa})$, for $\kappa\in W(\mu)$ with
$\mu\in\mathcal{P}^C_+$, is computed in the same way as in the
$B_N$-case.

\subsubsection{The case $D_N$}
The Weyl group is given by $W=\Sigma_N\ltimes \mathbb{Z}_2^{N-1}$,
and the action of $w=(\sigma,\varepsilon)\in W$ on a weight
$\lambda\in\mathcal{P}^+_D$ is given by Eq. \eqref{waction} with
$\varepsilon_j\in \{ 1 ,-1\}$ for $j=1,\ldots ,N$ such that
$\varepsilon_1\varepsilon_2\cdots \varepsilon_N=1$.

To construct the $D_N$-type Macdonald polynomial $p_\lambda^{D}$
associated to a weight $\lambda\in\mathcal{P}^D_+$, one employs
Theorem \ref{m-dc:thm} with
\begin{subequations}
\begin{eqnarray}
\epsilon_\mu^D &=& \sum_{j=1}^N \left (t^{2N-j-1}
q^{\mu_j} + t^{j-1} q^{-\mu_j} \right) , \\
\rho_D &= &(N-1,N-2,\ldots ,1,0).
\end{eqnarray}
\end{subequations}
For a weight $\kappa\in W(\mu)$ with $\mu\in\mathcal{P}^+_D$, the
sign $\det (w_{\rho+\kappa})$ is given by the signature of the
shortest permutation rearranging $\rho+\kappa$ such that the
absolute values of its parts become weakly decreasing.

\begin{remark}[i]
Example: for $R=D_3$ and $\lambda = (2,1,0)$ the determinantal
formula reads
\begin{small}
\begin{equation*}
\begin{split}
& \qquad \qquad
\makebox{\begin{normalsize}$p_{2,1,0}$\end{normalsize}}=
\bigl((\epsilon^D_{{2,1,0}}-\epsilon^D_{{1,1,1}})
(\epsilon^D_{{2,1,0}}-\epsilon^D_{{1,1,-1}})(\epsilon^D_{{2,1,0}}
-\epsilon^D_{{1,0,0}})\bigr)^{-1} \times \\
& \left |\begin {array}{cccc} m_{{1,0,0}}&\epsilon^D_{{2,1,0}}-
\epsilon^D_{{1,0,0}}&0&0
\\ \noalign{\medskip}m_{{1,1,-1}}&-\epsilon^D_{{2,1,0}}
+\epsilon^D_{{1,-1,1}}&\epsilon^D_{{2,1,0}}
-\epsilon^D_{{1,1,-1}}&0\\\noalign{\medskip}m_{{1,1,1}}&
-\epsilon^D_{{2,1,0}}+\epsilon^D_{{1,-1,- 1}}
&0&\epsilon^D_{{2,1,0}}-\epsilon^D_{{1,1,1}}\\\noalign{\medskip}
m_{{2,1,0}}&-\epsilon^D_{{1,-2,0}}+\epsilon^D_{{-1,2,0}}
&-2\,\epsilon^D_{{2,1,0}}+\epsilon^D_{{1,0,-2}}+
\epsilon^D_{{0,2,-1}}&-2\,\epsilon^D_{{2,1,0}}+\epsilon^D_{{1,0,2}}
+\epsilon^D_{{0,2,1}}\end {array}\right | ,
\end{split}
\end{equation*}
\end{small}
with
$\epsilon^D_{m_1,m_2,m_3}=(t^{4}q^{m_1}+q^{-m_1})
+(t^3q^{m_2}+t\,q^{-m_2})+(t^2q^{m_3}+t^2q^{-m_3})$.
\end{remark}

\begin{remark}[ii]
For the root systems $B_N$ and $C_N$ the minuscule weights $\pi$
of the dual root systems are unique. Specifically, for $R=B_N$ the
minuscule weight is given by lowest fundamental weight $\omega_1$
of $R^\vee$ $(=C_N$), and for $R=C_N$ it is given by the highest
fundamental weight $\omega_N$ of $R^\vee$ ($=B_N$). In the cases
of the root systems $A_{N-1}$ and $D_N$ the above formulas for the
eigenvalues $\epsilon_\mu$ correspond to picking for the minuscule
weight $\pi$ the lowest fundamental weight $\omega_1$ of $R^\vee$
($=R$). However, in these cases there exist actually several
alternative possibilities for the choice of the minuscule weight
$\pi$. For $R=A_{N-1}$ we could work with each of the eigenvalues
\begin{equation*}
\epsilon^A_{\mu ,r} = \sum_{\begin{subarray}{c} J\subset \{
1,\ldots , N\} \\ |J|=r \end{subarray}} \prod_{j\in J} t^{N-j}
q^{\lambda_j}, \;\;\;\; r=1,\ldots ,N-1,
\end{equation*}
corresponding to the fundamental weights $\omega_1,\ldots
,\omega_{N-1}$, respectively. For $R=D_N$ we could alternatively
work with the eigenvalues
\begin{equation*}
\epsilon^D_{\mu ,N-1}=(\epsilon^{D}_{\mu ,+ }- \epsilon^{D}_{\mu,
-})/2  \;\;\;\;\text{or} \;\;\;\;
\epsilon^D_{\mu,N}=(\epsilon^{D}_{\mu,+ }+\epsilon^{D}_{\mu,-})/2,
\end{equation*}
where
\begin{eqnarray*}
\epsilon^{D}_{\mu ,+ } &= & q^{-\frac{1}{2}\sum_{j=1}^N \mu_j}
\prod_{j=1}^N \left(
t^{N-j}q^{\mu_j}+1\right) ,\\
\epsilon^{D}_{\mu,-} &= & q^{-\frac{1}{2}\sum_{j=1}^N \mu_j}
\prod_{j=1}^N \left( t^{N-j}q^{\mu_j}-1\right)  ,
\end{eqnarray*}
corresponding to the fundamental (spin) weights $\omega_{N-1}$ and
$\omega_N$, respectively. From a computational point of view,
however, in this last case it is more efficient to work simply
with the eigenvalues $\epsilon^{D}_{\mu, + }$, corresponding to
the linear combination of Macdonald operators
$D=D_{\omega_{N-1}}+D_{\omega_N}$. For instance, in the example of
Remark (i) above, this amounts to replacing the eigenvalues
$\epsilon^D_{m_1,m_2,m_3}$ by
$\epsilon^D_{m_1,m_2,m_3,+}=q^{-(m_1+m_2+m_3)/2}(t^{2}\,q^{m_1}+1)
(t\,q^{m_2}+1)(q^{m_3}+1)$.
\end{remark}

\section{Macdonald Polynomials: the case of general $t_\alpha$}\label{sec7}
In this section we will briefly indicate how to generalize the
results of the previous section to the case of Macdonald
polynomials with general parameters $t_\alpha$ such that
$t_{w(\alpha)}=t_\alpha$ for all $w\in W$.  We will keep the
restriction that our root system $R$ is reduced and that the dual
root system $R^\vee$ has a minuscule weight $\pi$.

For general $W$-invariant $t_\alpha$-parameters the Macdonald
operator becomes \cite{mac:orthogonal}
\begin{equation}
D_\pi = \frac{1}{|W_\pi|} \sum_{w\in W} \Bigl( \prod_{\alpha\in
R^+} \frac{1-t_\alpha^{\langle \pi ,\alpha\rangle} e^{w(\alpha)}}
     {1-e^{w(\alpha)}} \Bigr) T_{w(\pi), q} .
\end{equation}
It is convenient to reparameterize the $t_\alpha$ as
\begin{equation*}
t_\alpha =q^{g_\alpha}
\end{equation*}
(with $g_{w(\alpha )}=g_\alpha$, $\forall w\in W$). The action of
$D_\pi$ on the monomial basis can be written as \cite[Section
5]{mac:orthogonal}
\begin{equation}\label{mac-action}
D_\pi\, m_\lambda = q^{\langle \pi ,\rho_g\rangle} \sum_{ X\subset
R^+} (-1)^{|X|} \sum_{\nu \in W(\lambda)}
 q^{\langle \pi ,\nu+\rho_g(X^c) -\rho_g (X)  \rangle}
\chi_{\nu-2\rho (X)} ,
\end{equation}
where $X^c=R^+\setminus X$ and
\begin{equation*}
\rho (X) = \frac{1}{2}\sum_{\alpha\in X} \alpha ,\;\;\;\; \;\;\;\;
\rho_g (X)=\frac{1}{2}\sum_{\alpha \in X} g_\alpha\, \alpha
\end{equation*}
(so $\rho =\rho (R^+)$ and $\rho_g=\rho_g(R^+)$). Bringing the
action in Eq. \eqref{mac-action} to triangular form gives
\begin{equation}
D_\pi\, m_\lambda = \epsilon_\lambda\, \chi_\lambda +
\sum_{\mu\in\mathcal{P}^+,\,\mu\prec\lambda} b_{\lambda \mu}\,
\chi_\mu ,
\end{equation}
with
\begin{eqnarray*}
\epsilon_\lambda &=&  q^{\langle \pi ,\rho_g\rangle}
\sum_{\tau\in W(\pi )} q^{\langle \tau , \lambda +\rho_g\rangle}, \\
b_{\lambda \mu}  &=& \\
&& \!\!\!\!\!\!\! \sum_{\begin{subarray}{c}
\nu\in W(\lambda) \cap (W(\mu +\rho )+\rho(X)-\rho(X^c)) \\
X\subset R^+
\end{subarray}}
(-1)^{|X|} \det (w_{\nu+\rho(X^c)-\rho(X)}) q^{\langle \pi ,\nu
+\rho_g(X^c)-\rho_g(X) \rangle}.
\end{eqnarray*}

We thus wind up with the following determinantal construction of
the Macdonald polynomials for general $W$-invariant parameters.
\begin{theorem}[Determinantal Construction]\label{mac-dc:thm}
For $\lambda\in\mathcal{P}^+$, let
\begin{equation*}
p_\lambda = m_\lambda +
\sum_{\mu\in\mathcal{P}^+,\,\mu\prec\lambda} c_{\lambda \mu}\,
m_\mu
\end{equation*}
denote the (monic) Macdonald polynomial with
$t_\alpha=q^{g_\alpha}$. Then we have---upon setting for $\mu
,\nu\in\mathcal{P}^+$
\begin{eqnarray*}
\epsilon_\mu &=&  q^{\langle \pi ,\rho_g\rangle}
\sum_{\tau\in W(\pi )} q^{\langle \tau , \mu +\rho_g\rangle},\\
d_{\mu\nu} &=& b_{\mu\nu}-\epsilon_\lambda a_{\mu\nu},
\end{eqnarray*}
with
\begin{eqnarray*}
a_{\mu\nu} \!\!&=&\!\! \sum_{\kappa\in W(\mu) \cap (W(\nu +\rho
)-\rho)}
\det (w_{\rho +\kappa}) ,\\
b_{\mu\nu }  \!\! &=&\!\! \\
&& \!\!\!\!\!\!\! \sum_{\begin{subarray}{c}
\kappa\in W(\mu) \cap (W(\nu +\rho )+\rho(X)-\rho(X^c)) \\
X\subset R^+
\end{subarray}}
(-1)^{|X|} \det (w_{\kappa+\rho(X^c)-\rho(X)}) q^{\langle \pi
,\kappa +\rho_g(X^c)-\rho_g(X) \rangle}
\end{eqnarray*}
(so $d_{\mu\mu}=\epsilon_\mu-\epsilon_\lambda$)---that:
\begin{itemize}
\item[i)] the polynomial $p_\lambda$ is represented explicitly
by the determinantal formula in Theorem \ref{df:thm},
\item[ii)] the coefficients $c_{\lambda \mu}$ of its monomial
expansion are generated by the linear recurrence in Corollary
\ref{lrr:cor},
\item[iii)] the expansion coefficients $c_{\lambda \mu}$ are
given in closed form by the formula in Corollary \ref{eme:cor}.
\end{itemize}
\end{theorem}

\begin{remark}[i]
{}From a computational standpoint the formulas of Theorem
\ref{mac-dc:thm} are much less effective than the determinantal
constructions for the $(q,t)$ Macdonald polynomials (Theorem
\ref{m-dc:thm}) and (especially) for the Heckman-Opdam polynomials
(Theorem \ref{ho-dc:thm}). This is because the action of the
general Macdonald operator on the monomial basis (cf. Eq.
\eqref{mac-action}) is much less sparse and the matrix elements
are moreover much more complex than in these two previous cases.
This renders Theorem \ref{mac-dc:thm} presumably only of limited
practical value.
\end{remark}

\begin{remark}[ii]
The most general class of Macdonald polynomials admits a richer
parameter structure connected with admissible pairs of root
systems $(R,S)$ \cite{mac:orthogonal}. (From this perspective the
polynomials studied here are of the type $(R,R)$.) Since Macdonald
in fact gives the action of the Macdonald operator on the monomial
basis for general admissible pairs, it is not difficult to
generalize Theorem \ref{mac-dc:thm} also to this context (at the
expense of having to introduce a more elaborate notational
apparatus).
\end{remark}

\appendix
\section{Determinants of Hessenberg Matrices}\label{appA}
In this appendix we recall a classic recursive method for the
efficient evaluation of the determinant of a Hessenberg matrix
\cite{wil:algebraic}.  This recursive method was used in Section
\ref{sec3} for the explicit evaluation of our determinantal
formula for the eigenbasis of the regular triangular operators in
$\mathcal{A}^W$.

\begin{lemma*}[\cite{wil:algebraic}]
Let $|D|$ be the Hessenberg determinant
 \begin{equation*}
  |D| \,=\,
  \begin{vmatrix}
  m_1&-d_{1,2} & 0 & \cdots & 0\\
  m_2 &d_{2,2}& -d_{2,3}& \ddots & 0 \\
  \vdots &\vdots &\vdots&\ddots & 0 \\
  \makebox[1ex]{}m_{n-1} & d_{n-1,2} &d_{n-1,3} &\cdots & -d_{n-1,n}\\
  m_n&d_{n,2} & d_{n,3}& \cdots &d_{n,n}
  \end{vmatrix},
  \end{equation*}
 with nonzero elements on the super-diagonal: $d_{j-1,j}\neq 0$ ($1<j\leq n$).
  Then the expansion of $|D|$ with respect to the first column is of the form
  \begin{equation*}
  |D|=\sum_{\ell=1}^n\, c_\ell\, m_\ell,
  \end{equation*}
  with
  $c_n=d_{1,2}d_{2,3}\cdots d_{n-1,n}$ and
  \begin{equation*}
  c_{\ell-1}=\frac{1}{d_{\ell-1,\ell}}\sum_{j=\ell}^n c_j\, d_{j,\ell}
  \end{equation*}
  ($1<\ell\leq n$).
  \end{lemma*}
\begin{proof} (\cite{lap-las-mor:determinantal2})
We denote the columns of our Hessenberg matrix by the
$n$-dimen\-sional vectors $\mathbf{m}$ and $\mathbf{d}^{(k)}$
($1<k\leq n$), respectively. Let $\mathbf{c}=(c_1,c_2,\ldots
,c_n)$ be a nonzero vector perpendicular to the hyperplane spanned
by the $(n-1)$ linear independent columns $\mathbf{d}^{(2)},\ldots
,\mathbf{d}^{(n)}$. Because the value of the determinant is equal
to the (hyper)volume of the polygon determined by
$\mathbf{d}^{(2)},\ldots ,\mathbf{d}^{(n)}$, multiplied by the
height of $\mathbf{m}$ in the perpendicular direction $\mathbf{c}$
(possibly up to a sign), we conclude that $|D|$ must be
proportional to the scalar product of $\mathbf{c}$ and the first
column $\mathbf{m}$, i.e. $|D|\sim
(\mathbf{c},\mathbf{m})=\sum_{\ell=1}^n c_\ell\, m_\ell$.
Furthermore, since the cofactor of $m_n$ in $|D|$ equals the
product of the elements on the super-diagonal, one sees that in
fact $|D|=\sum_{\ell=1}^n c_\ell\, m_\ell$ provided the
normalization of $\mathbf{c}$ is fixed such that
$c_n=d_{1,2}d_{2,3}\cdots d_{n-1,n}$. The lemma now follows from
the observation that the requirement that $\mathbf{c}$ be
orthogonal to the columns $\mathbf{d}^{(2)},\ldots
,\mathbf{d}^{(n)}$ translates itself directly into the stated
recurrence relations for the components $c_j$, $j=1,\ldots ,n$.
\end{proof}

\section{Counting the Orders of Stabilizer Subgroups}\label{appB}
To build the determinantal formula for the Heckman-Opdam
polynomials in Theorem \ref{ho-dc:thm}, one frequently needs to
compute the orders of stabilizer subgroups of the Weyl group. In
this appendix we include a short proof of a useful formula for the
orders of these stabilizers that can be found e.g. in Ref.
\cite[Section 12]{mac:orthogonal}.

\begin{proposition}[Orbit Size]\label{orbit:prp}
Let $\lambda\in \mathcal{P}$. Then the size of the Weyl orbit
through $\lambda$ is given by
\begin{equation*}
|W(\lambda )| =
\prod_{\begin{subarray}{c} \alpha\in R^+ \\
             \langle \lambda ,\alpha^\vee\rangle \neq 0
       \end{subarray}}
\frac{\langle \rho ,\alpha^\vee\rangle
+1+\frac{1}{2}\delta_{\frac{\alpha}{2}}} {\langle \rho
,\alpha^\vee\rangle +\frac{1}{2}\delta_{\frac{\alpha}{2}}} ,
\end{equation*}
where $\rho =\frac{1}{2}\sum_{\alpha\in R^+}\alpha$, and
\begin{equation*}
\delta_\alpha =
\begin{cases}
1 &\text{if}\;\; \alpha \in R ,\\
0& \text{otherwise} .
\end{cases}
\end{equation*}
\end{proposition}
\begin{proof}
Let us first assume that $\lambda$ is dominant. The orbit size can
be obtained in this case from the following evaluation formula for
the Heckman-Opdam polynomial $p_\lambda$ at $x=0$ (i.e. at the
identity element of the torus $\mathbb{T}=E/(2\pi
\mathcal{Q}^\vee)$) \cite{hec-sch:harmonic}
\begin{equation*}
p_\lambda (0) = \prod_{\alpha\in R^+} \frac{[\langle \rho_g
,\alpha^\vee\rangle +g_\alpha
            + \frac{1}{2} g_{\frac{\alpha}{2}}]_{\langle \lambda ,\alpha^\vee\rangle}}
     {[\langle \rho_g ,\alpha^\vee\rangle
            + \frac{1}{2} g_{\frac{\alpha}{2}}]_{\langle \lambda ,\alpha^\vee\rangle}} ,
\end{equation*}
where $[a]_m=a(a+1)\cdots (a+m-1)$ (with the convention that
$[a]_0=1$), and by definition $g_{\alpha}=0$ if $\alpha\not\in R$.
Indeed, setting $g_\alpha=g$ for all $\alpha\in R$, and performing
the limit $g\to 0$, readily entails the formula of the
proposition. (Here one uses that $\lim_{g_\alpha\to 0} p_\lambda =
m_\lambda$ and that $m_\lambda (0)= |W(\lambda )|$.) The extension
to non-dominant weights $\lambda$ is immediate (cf. also the two
corollaries below).
\end{proof}

By picking $\lambda$ regular (for instance strongly dominant), one
gets a formula for the order of the Weyl group.

\begin{corollary}[Order of the Weyl Group]\label{weyl:cor}
The order of the Weyl group is given by
\begin{equation*}
|W | = \prod_{\alpha\in R^+} \frac{\langle \rho
,\alpha^\vee\rangle +1+\frac{1}{2}\delta_{\frac{\alpha}{2}}}
{\langle \rho ,\alpha^\vee\rangle
+\frac{1}{2}\delta_{\frac{\alpha}{2}}} .
\end{equation*}
\end{corollary}

Dividing the order of the Weyl group by the size of the orbit
produces the order of the stabilizer subgroup.
\begin{corollary}[Order of the Stabilizer]\label{stabilizer:cor}
Let $\lambda\in\mathcal{P}$. The order of the stabilizer of
$\lambda$ is given by
\begin{equation*}
|W_\lambda | =
\prod_{\begin{subarray}{c} \alpha\in R^+ \\
             \langle \lambda ,\alpha^\vee\rangle = 0
       \end{subarray}}
\frac{\langle \rho ,\alpha^\vee\rangle
+1+\frac{1}{2}\delta_{\frac{\alpha}{2}}} {\langle \rho
,\alpha^\vee\rangle +\frac{1}{2}\delta_{\frac{\alpha}{2}}} .
\end{equation*}

\begin{remark}[i]
It is clear that the formulas of Proposition \ref{orbit:prp} and
Corollary \ref{stabilizer:cor} in fact serve to compute the sizes
of the orbit $W(x)$ and the stabilizer $W_x$ for {\em any} vector
$x\in E$ (i.e. not just weight vectors).
\end{remark}

\begin{remark}[ii]
For a reduced root system $R$ the above formulas simplify somewhat
as $\delta_{\frac{\alpha}{2}}=0$ in this case.
\end{remark}

\end{corollary}

\bibliographystyle{amsplain}

\end{document}